\documentclass[11pt,bull-l]{amsart}
\usepackage[top = 0.9in, bottom = 0.9in, left =1in, right = 1in]{geometry}
\usepackage{geometry}

\usepackage{comment}
\usepackage{url}
\usepackage{tikz}
\usetikzlibrary{shapes.geometric, arrows}
\usepackage{lipsum}
\usepackage{amsfonts,amscd,amssymb,amsmath,mathrsfs}
\usepackage{graphicx,cite}
\usepackage{epstopdf}
\usepackage{tabulary}
\usepackage{booktabs}
\usepackage{tikz}
\usepackage{pgfplots}
\usepackage{wrapfig}
\usepackage[font = footnotesize]{caption}
\usepackage{wrapfig}
\usepackage{cutwin}
\usepackage{algpseudocode}
\usepackage{algorithm}
\usepackage{subcaption}
\ifpdf
  \DeclareGraphicsExtensions{.eps,.pdf,.png,.jpg}
\else
  \DeclareGraphicsExtensions{.eps}
\fi

\tikzstyle{startstop} = [rectangle, rounded corners, minimum width=3cm, minimum height=1cm,text centered, draw=black, fill=red!30]
\tikzstyle{io} = [trapezium, trapezium left angle=70, trapezium right angle=110, minimum width=3cm, minimum height=1cm, text centered, draw=black, fill=blue!30]
\tikzstyle{process} = [rectangle, minimum width=3cm, minimum height=1cm, text centered, draw=black, fill=orange!30]
\tikzstyle{decision} = [diamond, minimum width=3cm, minimum height=1cm, text centered, draw=black, fill=green!30]
\tikzstyle{arrow} = [thick,->,>=stealth]

\newtheorem{theorem}{Theorem}[section]

\newtheorem{remark}{Remark}[section]

\DeclareMathOperator{\imag}{Im}

\usepackage{enumitem}
\setlist[enumerate]{leftmargin=.5in}
\setlist[itemize]{leftmargin=.5in}

\title{An artificially-damped Fourier method for dispersive evolution equations}

\author{Anne Liu}
\address{University of Washington, Seattle, WA}
\email{aliu7@uw.edu}
  
\author{Thomas Trogdon}
\address{University of Washington, Seattle, WA}
\email{trogdon@uw.edu}

\thanks{TT is partially supported by NSF DMS-1945652.}    

\usepackage{amsopn}

\let\vec=\mathbf
\newcommand{\ii}{\mathrm i}
\newcommand{\ee}{\mathrm e}
\newcommand{\dd}{\mathrm d}

\usepackage{hyperref}
\begin{document}

\maketitle

% REQUIRED
\begin{abstract}
Computing solutions to partial differential equations using the fast Fourier transform can lead to unwanted oscillatory behavior. Due to the periodic nature of the discrete Fourier transform, waves that leave the computational domain on one side reappear on the other and for dispersive equations these are typically high-velocity, high-frequency waves. However, the fast Fourier transform is a very efficient numerical tool and it is important to find a way to damp these oscillations so that this transform can still be used. In this paper, we accurately model solutions to four nonlinear partial differential equations on an infinite domain by considering a finite interval and implementing two damping methods outside of that interval: one that solves the heat equation and one that simulates rapid exponential decay. Heat equation-based damping is best suited for small-amplitude, high-frequency oscillations while exponential decay is used to damp traveling waves and high-amplitude oscillations. We demonstrate significant improvements in the runtime of well-studied numerical methods when adding in the damping method.
\end{abstract}

\section{Introduction}

Distortions to numerical solutions of partial differential equations (PDEs) on the line caused by employing the periodic Fourier method are often mitigated by increasing the size of the computational domain. While this does mitigate the errors caused by the periodic copies generated by the discrete Fourier transform, it also increases the runtime. In this paper, we outline a method of artificial damping that allows for a smaller computational domain and significantly shorter runtimes. We demonstrate the effectiveness and generality of the method by showing how it applies to the Korteweg-de Vries (KdV) equation, the nonlinear Schrödinger (NLS) equation, a Riemann problem for the KdV equation, and a Riemann problem for the Kawahara equation\footnote{All computations in this paper, with the exception of computing reference solutions for the KdV, NLS and Kawahara equations, were done on a Microsoft Surface Pro 6 Laptop with 8GB RAM, 256GB SSD, a 1.60GHz CPU, the Intel® Core™ i5-8250U Processor, and the Windows 10 operating system.}. The Jupyter notebooks used for computations in this paper, as well as a template for general implementation of this damping method, are provided in a GitHub repository \cite{Liu_Artificially-damped_Fourier_method}.

We also present a heuristic error analysis of the method as applied to the linearized KdV equation. Artificially damping of solutions of PDEs has been discussed in other contexts and referred to by many other names. In \cite{sponge}, Carmigniani and Violeau model solutions to water wave problems in an open ocean by adding terms that introduce linear dissipative forces according to a sponge layer function. Unlike our method, which is intended for general dispersive evolution equations, this method is intended for full Navier-Stokes solvers. In \cite{beach}, Cao, Beck, and Schultz discuss the modeling of free surface waves on an infinite domain by adding a damping term to the free surface boundary condition. They refer to the region of damping as an absorbing beach. In this work numerical results for a damping term added to Bernoulli's equation are given, but the discussion of the method's generality and accuracy for nonlinear PDEs is limited. 

For our method, we use the fast Fourier transform \cite{ogft}, described in Section~\ref{sec:fourier}, and the standard fourth-order Runge-Kutta method \cite{leveque} to evolve the Fourier coefficients forward in time. Fourier series are inherently periodic, so when the Fourier transform of an aperiodic function on \(\mathbb{R}\) is computed by restricting it to an interval\footnote{A note on intervals: \([-L,L]\) is the interval on which the solution is computed, \([-R,R]\) is the interval on which we will compare the solution to the true solution, and \([-L,-P_-]\) and \([P_+,L]\) are the damping regions at the edges of the computational domain. The sizes of the intervals are related as follows: \(-L<-P_-<-R\) and \(R<P_+<L\).} \([-L, L]\), the inverse transform on the whole real line becomes periodic with a period of \(2L\). This is illustrated in Figure~\ref{fig:fig2}. This forced periodicity distorts solutions and becomes increasingly problematic at larger times. Typically, a large enough interval is used to minimize interference from the periodic copies of the solution outside the interval. However, depending on the problem, a prohibitively large interval may be required. For this reason, we seek to artificially damp behavior outside the interval of interest, \([-R,R]\). Damping the solution on \([-L,-P_-]\) and \([P_+,L]\) creates a buffer between periodic solutions and stops neighboring solutions from interfering with the desired solution on \([-R,R]\).

We have two primary methods of artificial damping. The first is to solve the heat equation using Strang-splitting\footnote{Lower-order or higher-order splitting methods could be considered but we use Strang-splitting due to its wide use.} outside \([-R,R]\). This technique gradually damps oscillations and is best suited for small-amplitude, high-frequency oscillations. Our second method is to add exponential decay to solutions outside \([-R,R]\). This technique is more aggressive and is best suited for solitons, traveling waves, and high-amplitude oscillations. Both methods are discussed in more detail in Section~\ref{sec:damping1}.

In Section~\ref{sec:KdV}, we solve the KdV equation \cite{ogkdv}, which has a number of physical applications, including shallow-water waves \cite{shallow}, internal solitary waves in the ocean \cite{internal}, and ion acoustic waves \cite{ion}. It is the first known PDE that is solvable with the inverse scattering transform \cite{ist} and has also been solved using numerical inverse scattering without any boundary approximation \cite{trogdonschrodinger} and we use this method to compute a reference, or benchmark, solution. We consider a solution with the initial condition, \(q(x,0)=1.3\ee^{-x^2}\), that produces a dispersive tail and a single soliton. We numerically compute solutions without the damping technique on intervals of varying sizes to demonstrate the pitfalls of the periodic Fourier method. We then artificially damp the solution outside the interval of interest and are able to compute an accurate solution on a smaller computational domain, thus significantly decreasing the runtime. This section also includes a discussion of how the accuracy and runtime of the method depends on the damping parameters, as well as some error bound analysis on the damping method for the linearized KdV equation \eqref{eq:linkdv}.

In Section~\ref{sec:NLS}, we solve the NLS equation \cite{ognls}, which can be used to model Bose-Einstein condensation and plasma physics \cite{bose}, deep water waves \cite{deep}, and nonlinear optics \cite{optics}. There exist many numerical methods for solving the NLS equation, including numerical Fourier methods \cite{split}, a relaxation scheme \cite{relax}, and numerical inverse scattering \cite{trogdonschrodinger} which is again used to compute a reference solution. We compute solutions with the initial condition, \(q(x,0)=(1+x)\ee^{\ii x-0.7x^2}\), which produces solitons traveling in both directions and high-amplitude oscillations. 

\begin{remark}
We consider the integrable KdV and NLS equations because of the existence of the numerical inverse scattering transform that enables one to compute high-accuracy reference solutions that do not suffer from boundary approximation issues.  This is analogous to numerically evaluating the Fourier transform solution for a linear PDE on the line.  Thus, one can estimate the true error of a numerical method applied to these equations to a high degree of accuracy.  This approach first appeared in \cite{Bilman2017a}.
\end{remark}

In Section~\ref{sec:shock}, we compute solutions to a Riemann problem for the KdV equation, motivated by a problem discussed by, for example, Grava and Klein in \cite{cauchyproblem}. Solutions to this problem are dispersive shock waves, which have been observed in plasmas \cite{shockplasma}, optical fibers \cite{ballistic}, and viscously deformable media \cite{viscous}. We use an aperiodic step-like function for the initial condition, \(q(x,0)=\frac{1}{1+\ee^{10x}}\), which has a nearly periodic derivative. Inspired by \cite{Sprenger2017, hoefershearer}, we modify the PDE to instead solve for the derivative of the desired solution. 

Similarly, in Section~\ref{sec:kawa}, we solve the Kawahara equation with an aperiodic step-like initial condition, \(q(x,0)=\frac{1}{1+\ee^{10x}}-1\), also producing a dispersive shock wave. Therefore, we utilize the approximately periodic derivative of the initial condition in the same way as with the Riemann problem for the KdV equation. The problem we choose to solve is motivated by work done by Sprenger and Hoefer in \cite{Sprenger2017}. Physical applications of the Kawahara equation include magneto-acoustic waves in plasmas \cite{ogkawahara} and shallow water waves with surface tension \cite{tension}. By solving these four nonlinear dispersive PDEs, we demonstrate the effectiveness and generality of the damping technique, which we believe can be applied to a wide class of dispersive PDEs with the appropriate modifications.

\section{The Fourier method}
\label{sec:fourier}

We begin with an elementary overview of the classical method on which we expand. We consider a \(k\)th-order partial differential equation for \(q(x,t)\) of the form
\begin{equation}
    q_t + \mathcal{L}q + \mathcal{N}(q,q_x) = 0,
    \label{eq:quasigeneral}
\end{equation}
where \(\mathcal{L} = \alpha_1 \partial_x + \alpha_2\partial_x^2 + \cdots + \alpha_k\partial_x^k\) is a linear combination of the spatial derivatives and \(\mathcal{N}\) is a nonlinear function\footnote{In effect, we require this function to be a polynomial.  This term could depend on higher derivatives of $q$ but this may prevent one from using explicit time stepping, complicating matters significantly.} of \(q\) and \(q_x\). We assume that if \(f(x) = e^{ikx}\), then \(\mathcal{L}f(x) = \omega(k)f(x)\), where \(\omega: \mathbb{R}\rightarrow i\mathbb{R}\) is purely imaginary. We use the Fourier method \cite{ogft} with approximations\footnote{In principle, one should use a power of $2$ for the number of Fourier modes to improve the efficiency of the FFT.  We will use $2J+1$ modes to simplify the presentation.}
\begin{align*}
    q(x,t) &\approx\mathcal{F}^{-1}_J(\vec c(t))(x):=\sum_{j=-J}^J c_j(t)\ee^{\frac{\ii j\pi}{L}x},\\
    \vec c(t) &=(c_{-J}(t),c_{-J+1}(t),\cdots, c_J(t)),
\end{align*}
and
\[c_j(t)\approx\frac{1}{2L}\int_{-L}^{L}q(x,t)\ee^{-\frac{\ii j\pi}{L}x}\dd x.\]
Here, \(\mathcal{F}_J\) is the discrete Fourier transform for functions on \([-L,L]\). The approximate coefficients of the initial condition, \(c_j(0)\), found using the FFT, are mathematically equivalent to using the trapezoidal rule \cite{traprule} and discretizing the computational domain, \([-L,L]\), using \(m = 2J + 1\) evenly-spaced grid points as follows:
\begin{equation}
    c_j(0) = \frac{1}{2J + 1}\sum_{i = 1}^{m}q_j(x_i,0)\ee^{-\frac{\ii j\pi}{L}x_i},
    \label{eq:trap}
\end{equation}
where
\begin{equation}
    x_i = -L + 2L\left(\frac{i-1}{m}\right).
    \label{eq:xi}
\end{equation}
We use the following finite-dimensional system of ordinary differential equations to approximate the Fourier coefficients of \(q(x,t)\):
\[\vec c'(t) + M \vec c(t) = F\left(\vec c(t)\right),\]
where 
\begin{align*}
    F(\vec c) &= -\mathcal{F}_J\left(\mathcal{N}\left(\mathcal{F}_J^{-1}(\vec c), \mathcal{F}_J^{-1}(D_J\vec c)\right)\right),\\
    M &= \alpha_1D_J + \cdots + \alpha_nD_J^k, \quad 
    D_J = \frac{\ii\pi}{L}\textrm{diag}(-J,-J +1,\cdots, J).
\end{align*}
Here we use the notation that a function $g : \mathbb R^n \to \mathbb R$ extends to a function
\begin{align*}
g: \underbrace{\mathbb R^m \times \cdots \times \mathbb R^m}_{n~\text{times}} \to \mathbb R^m,
\end{align*}
by applying $g$ entrywise
\begin{align*}
g(\vec x^{(1)}, \vec x^{(2)}, \ldots, \vec x^{(n)}) = \begin{bmatrix} g( x^{(1)}_1, x^{(2)}_1, \ldots, x^{(n)}_1) \\
\vdots \\
g( x^{(1)}_m, x^{(2)}_m, \ldots, x^{(n)}_m)\end{bmatrix}.
\end{align*}
This convention will be used throughout what follows. We then rewrite the system as
\[\ee^{-Mt}\left(\frac{\dd}{\dd t}\left(\ee^{Mt}\vec c(t)\right)\right) = F\left(\vec c(t)\right),\]
and define \(\vec a(t) = \ee^{Mt}\vec c(t)\) so that
\begin{equation}
    \vec a'(t) = \ee^{Mt}F\left(\ee^{-Mt}\vec a(t)\right).
    \label{eq:genode}
\end{equation}
Note that our assumptions on \(\mathcal{L}\) imply that \(M\) is diagonal with purely imaginary diagonal entries.
To compute solutions, we obtain the Fourier coefficients, \(\vec c(0)\), of the initial condition \(q(x,0)\) and advance them forward in time according to \eqref{eq:genode}. We use the standard RK4 method \cite{leveque}, as outlined in Algorithm~\ref{alg:rk4}. The algorithm takes in a vector \(\vec a^{(n)}\) that approximates \(\vec a(t_n)\) and computes the vector \(\vec a^{(n+1)}\) to approximate \(\vec a(t_{n+1}) = \vec a(t_n+\triangle t)\), where \(\triangle t\) is the size of one time step. Once we have obtained the approximate Fourier coefficients of the solution at the desired time, we use the inverse Fourier transform to compute the approximation of \(q(x,t)\) on the grid $(x_i)_{i=1}^m$. 
\begin{algorithm}
\caption{Fourth-order Runge-Kutta}
\label{alg:rk4}
\begin{algorithmic}[1]
\Function{\(\tt{rk4}\)}{$\vec a^{(n)},\;t_n,\;\triangle t$}
\State{$\vec f_1 = \ee^{Mt}F(\ee^{-Mt_n}\vec a^{(n)})$}
\State{$\vec f_2 = \ee^{M(t_n+\frac{\triangle t}{2})}F(\ee^{-M(t_n+\frac{\triangle t}{2})}(\vec a^{(n)}+\frac{\triangle t}{2}\vec{f_1}))$}
\State{$\vec f_3 = \ee^{M(t_n+\frac{\triangle t}{2})}F(\ee^{-M(t_n+\frac{\triangle t}{2})}(\vec a^{(n)}+\frac{\triangle t}{2}\vec{f_2}))$}
\State{$\vec f_4 = \ee^{M(t_n+\triangle t)}F(\ee^{-M(t_n+\triangle t)}(\vec a^{(n)}+\triangle t\vec{f_3}))$}
\State{$\vec a^{(n+1)} = \vec a^{(n)} + \frac{\triangle t}{6}(\vec f_1+2\vec f_2+2\vec f_3+\vec f_4)$}
\State{\textbf{return} $\vec a^{(n+1)}$}
\EndFunction
\end{algorithmic}
\end{algorithm}

\section{Artificial damping}
\label{sec:damping1}

Due to the periodic nature of the Fourier method, the approximate solution will have a period of \(2L\). As time evolves, non-zero behavior from neighboring solutions will enter the computational domain, \([-L, L]\), and distort the solution. Errors worsen over time as more collisions occur at the edges of the interval --- a result of the dispersive nature of the PDEs. This phenomenon is discussed in more detail in Section~\ref{sec:KdV}. In order to mitigate this problem, we introduce damping at the edges of the interval of computation\footnote{For all computations in this paper, \(P_-=P_+\), but we leave this distinction for generality.}: \([-L,-P_-]\) and \([P_+,L]\).
We modify the PDE in \eqref{eq:quasigeneral} to introduce damping:
\begin{equation}
    q_t + \mathcal{L}q + \mathcal{N}(q,q_x) = k_1\left(\sigma(x)q_{x}\right)_x-k_2(1-\gamma(x))q.
    \label{eq:modgeneral}
\end{equation}
The terms added to the right-hand side are effectively zero outside of the damping regions:
\[ k_1\left(\sigma(x)q_{x}\right)_x-k_2(1-\gamma(x))q \approx 0, \;\;\;\;\;|x\mp P_{\pm}|\gg 1.\]
We solve \eqref{eq:modgeneral} using splitting methods, as described in \cite{strang_split}. This requires two additional components. (1) We solve the heat equation,
\begin{equation*}
    q_t=k_1\left(\sigma(x)q_{x}\right)_x,
\end{equation*}
with a diffusion coefficient, \(k_1\sigma(x)\), that is chosen to be non-zero where damping is desired and zero elsewhere. Using the Fourier transform, we obtain an approximate ODE system for the Fourier coefficients of \(q(x,t)\),
\[\vec{c}'(t) = k_1D_J\mathcal{F}_J\left(\sigma(\vec{x})\mathcal{F}_J^{-1}\left(D_J\vec{c}(t)\right)\right).\]
Here, \(\vec{x}=(x_1,\cdots,x_m)\) is the vector of \(m\) evenly-spaced grid points on \([-L,L]\) as in \eqref{eq:xi}. We discretize the system using the trapezoidal method \cite{leveque}. This yields the following system of equations for \(\vec{c}^{(n)}\approx \vec{c}(t_n)\),
\begin{equation}
    B\vec{c}^{(n+1)} = A\vec{c}^{(n)},
    \label{eq:disc}
\end{equation}
where \(B\) and \(A\) are matrices defined by
\[B\vec{x} = \vec{x} - \frac{k_1\triangle t}{2}D_J\mathcal{F}_J\left(\Sigma \mathcal{F}^{-1}_N\left(D_J\vec{x}\right)\right)\]
and 
\[
A\vec{x} = \vec{x} + \frac{k_1\triangle t}{2}D_J\mathcal{F}_J\left(\Sigma \mathcal{F}^{-1}_N\left(D_J\vec{x}\right)\right),\]
with \(\Sigma=\textrm{diag}(\sigma(\vec{x}))\). We solve \eqref{eq:disc} using the conjugate gradient algorithm \cite{cg}, performing this diffusion every \(f_1\) time steps. We review the conjugate gradient algorithm in Algorithm~\ref{alg:cg}, where \(\vec{a} = A\vec{c}^{(n)}\) and \(\epsilon\) is a parameter that determines the accuracy of the approximation. We do not use a preconditioner or an initial guess in Algorithm~\ref{alg:cg} and the runtime of solutions could be improved through these additions. 
\begin{algorithm}
\caption{Conjugate gradient algorithm}
\label{alg:cg}
\begin{algorithmic}[1]
\Function{\tt{cg}}{$B,\; \vec{a},\;\epsilon$}
\State{$\vec{c}^{(n+1)} = \vec{0}; \; \vec{r} = \vec{a}; \; \vec{p} = \vec{r}; \; n = 0;$}
\While{$\textrm{norm}(\vec{r}) > \epsilon$}
    \State{$\vec{q} = B(\vec{p})$}
    \State{$\vec{a} = \frac{\vec{r}^*\vec{r}}{\vec{p}^*\vec{q}}$}
    \State{$\vec{c}^{(n+1)}$ += $\vec{a}\vec{p}$}
    \State{$\Bar{\vec{r}} = \vec{r}$}
    \State{$\vec{r}$ += $-\vec{a}\vec{q}$}
    \State{$\vec{a} = \frac{\vec{r}^*\vec{r}}{\Bar{\vec{r}}^*\Bar{\vec{r}}}$}
    \State{$\vec{p} = \vec{r} + \vec{a}\vec{p}$}
    \State{$n$ += $1$}
\EndWhile
\State{\textbf{return} $\vec{c}^{(n+1)}$}
\EndFunction
\end{algorithmic}
\end{algorithm}\\
(2) The second component of damping comes from solving 
\begin{equation*}
    q_t = k_2(1-\gamma(x))q.
\end{equation*}
We simulate rapid exponential decay as \(k_2\rightarrow\infty\) by judiciously multiplying the solution values by \(\gamma(x)\) every \(f_2\) time steps, where again \(\gamma(x)\) is chosen to enforce damping in the desired regions. This procedure is outlined in Algorithm~\ref{alg:pseudo}, where we define $\vec{q}^{(0)}\in\mathbb{R}^m$ to be a vector of \(q(x,0)\) evaluated on \(\vec{x} = (x_1,\cdots,x_m)\). Note that when calling \(\tt{rk4}\) in the algorithm, we use \(\tt{rk4}(\vec{c}^{(n)}, 0, \triangle t)\), as opposed to \(\tt{rk4}(\vec{a}^{(n)}, t, \triangle t)\). The PDEs we are solving are autonomous, so evolving the solution from \(t\) to \(t + \triangle t\) gives the same result as evolving it from 0 to \(\triangle t\). The latter is simpler computationally.

\begin{algorithm}
\caption{Artificial damping}
\label{alg:pseudo}
\textbf{Input:} A vector of the function values of the initial condition, $\vec{q}^{(0)}$. \\
\textbf{Output:} A vector of the function values at time \(T\), $\vec{q}^{(N)}$.

\begin{algorithmic}[1]
\State{$\vec{c}^{(0)} = \mathcal{F}_J(\vec{q}^{(0)})$} 
\State{$N = \frac{T}{\triangle t}$}
\For{\texttt{(n = 0; n < N; n++)}}
    \If{$\frac{n}{f_1}$ is an integer} \Comment{$\sigma(x)$ damping}
        \State{$\vec{a}^{(n+\frac{1}{2})} = \tt{rk4}(\vec{c}^{(n)}, 0, \frac{\triangle t}{2})$}
        \State{$\vec{\Tilde{c}} = \ee^{-M\frac{\triangle t}{2}}\vec{a}^{(n+\frac{1}{2})}$} \Comment{Strang-splitting}
        \State{$\vec{\Bar{c}} = \tt{cg}(B,A(\vec{\Tilde{c}}),\epsilon)$}
        \State{$\vec{a}^{(n+1)} = \tt{rk4}(\vec{\Bar{c}}, 0, \frac{\triangle t}{2})$}
        \State{$\vec{c}^{(n+1)} = \ee^{-M\frac{\triangle t}{2}}\vec{a}^{(n+1)}$}
    \Else \Comment{No damping}
        \State{$\vec{a}^{(n+1)} = \tt{rk4}(\vec{c}^{(n)}, 0, \triangle t)$}
        \State{$\vec{c}^{(n+1)} = \ee^{-M\triangle t}\vec{a}^{(n+1)}$}
    \EndIf
    \If{$\frac{n}{f_2}$ is an integer} \Comment{$\gamma(x)$ damping}
        \State{$\vec{c}^{(n+1)} = \mathcal{F}_J\left(\gamma(\vec{x})\cdot \mathcal{F}_J^{-1}(\vec{c}^{(n+1)})\right)$}
    \EndIf
\EndFor
\State{$\vec{q}^{(N)} = \mathcal{F}_J(\vec{c}^{(n)})$}\\
\Return $\vec{q}^{(N)}$
\end{algorithmic}
\end{algorithm}

\section{Korteweg-de Vries equation}
\label{sec:KdV}

\begin{figure}[htbp]
  \centering
  \includegraphics[width = 0.65\linewidth]{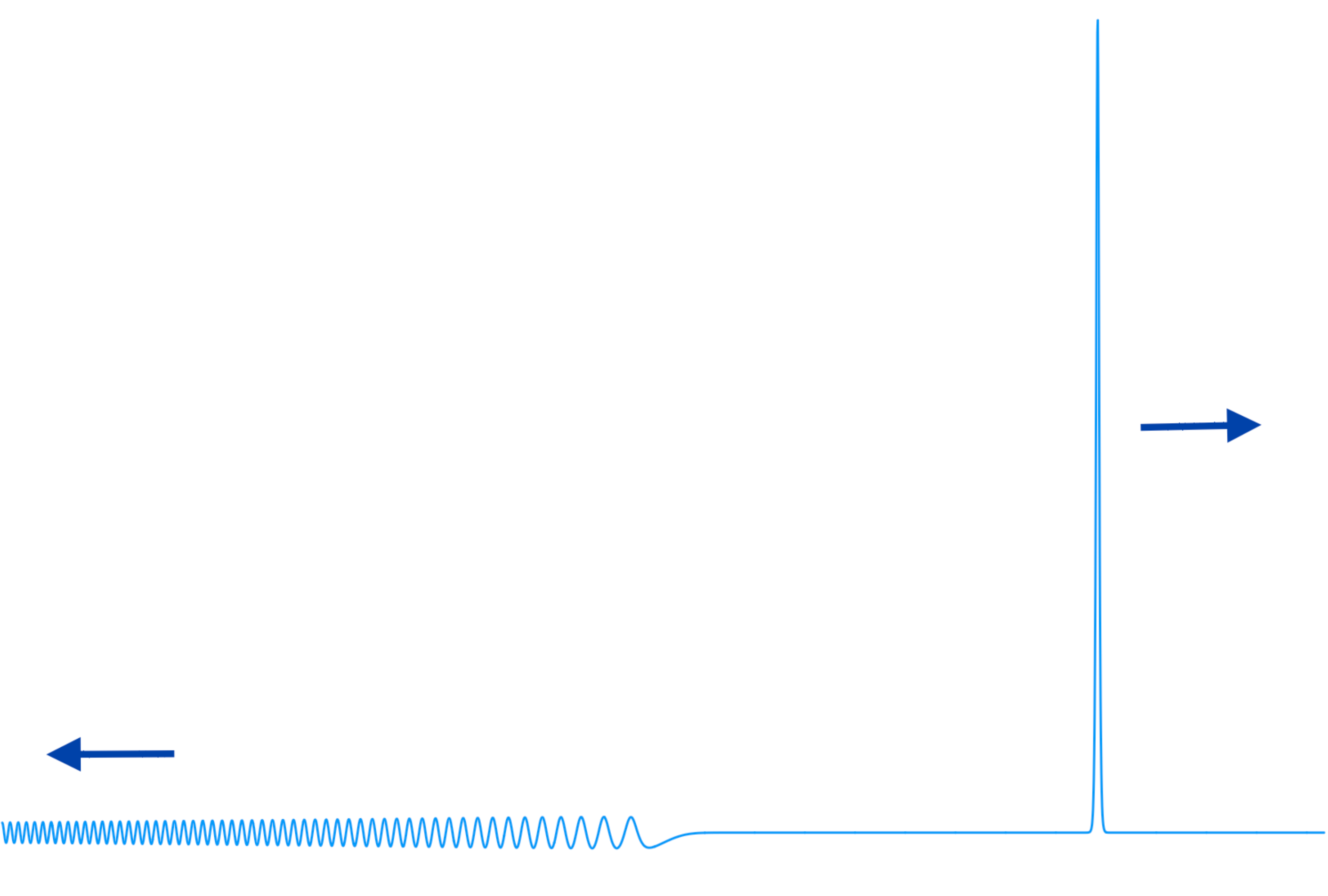}
  \caption{Solution to the KdV equation with initial condition \(q(x,0)=1.3e^{-x^2}\) at \(t = 150\).}
  \label{fig:fig1}
\end{figure}

We now apply our damping technique to the Korteweg-de Vries (KdV) equation \cite{ogkdv},
\begin{equation}
    q_t+6qq_x+q_{xxx}=0,
    \label{eq:kdv}
\end{equation}
with the initial condition \(q(x,0)=1.3\ee^{-x^2}\). This produces leftward traveling oscillations (i.e., a dispersive tail) and a rightward traveling soliton, as shown in Figure~\ref{fig:fig1}.

\subsection{The Fourier method for the KdV equation}
\label{sec:fourierkdv}

\begin{figure}[htbp]
  \centering
  \includegraphics[width = 0.8\linewidth]{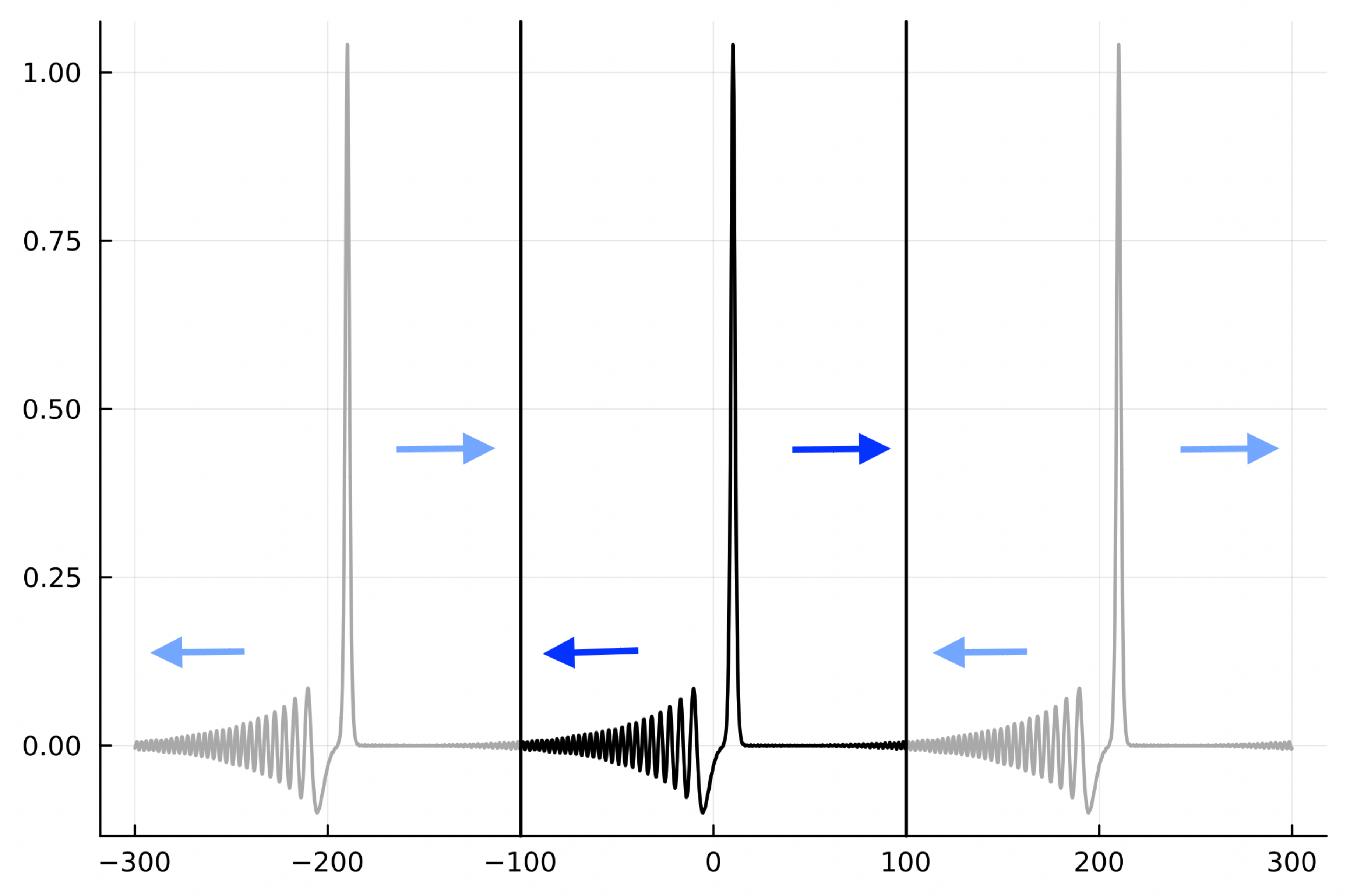}
  \caption{Solution to the KdV equation computed on \([-100,100]\) with initial condition \(q(x,0)=1.3e^{-x^2}\) at time \(t = 5\).}
  \label{fig:fig2}
\end{figure}

We seek to apply the Fourier method described in Section~\ref{sec:fourier}. The KdV equation is a third-order nonlinear PDE, so the general functions in \eqref{eq:quasigeneral} become
\(\mathcal{L}(q) = q_{xxx}\)
and
\(\mathcal{N}(q,q_x) = 6qq_x.\)
Therefore, \(M = D_J^3\) and the system of ordinary differential equations in \eqref{eq:genode} becomes
\[\vec{a}'(t) = \ee^{D_J^3t}F(\ee^{-D_J^3t}\vec{a}(t))\]
where \(\vec{a}(t) = \ee^{D_J^3t}\vec{c}(t)\),
\[F(\vec{c}) = -6\mathcal{F}_J\left(\mathcal{F}_J^{-1}(\vec{c}(t))\cdot\mathcal{F}_J^{-1}(D_J \vec{c}(t))\right),\]
and $\cdot$ denotes entrywise multiplication of vectors.
As was briefly discussed in Section~\ref{sec:fourier}, the Fourier method forces the solution to be periodic on the whole real line. In Figure~\ref{fig:fig2}, we plot the inverse transform of the Fourier coefficients of a solution to the KdV equation computed on \([-100,100]\). Outside of the computational domain, the solution has become periodic. 
At the bounds of the interval, the solution collides with the neighboring periodic solutions. The dispersive tail of the solution on \([100,300]\) flows into the interval from the right and introduces oscillations that are not present in the true solution. At a later time, the effective soliton in the solution on \([-300,-100]\) will enter the interval from the left and cause large errors.

\subsection{Undamped solutions of the KdV equation}

\begin{figure}[htbp]
  \centering
  \includegraphics[width = 0.8\linewidth]{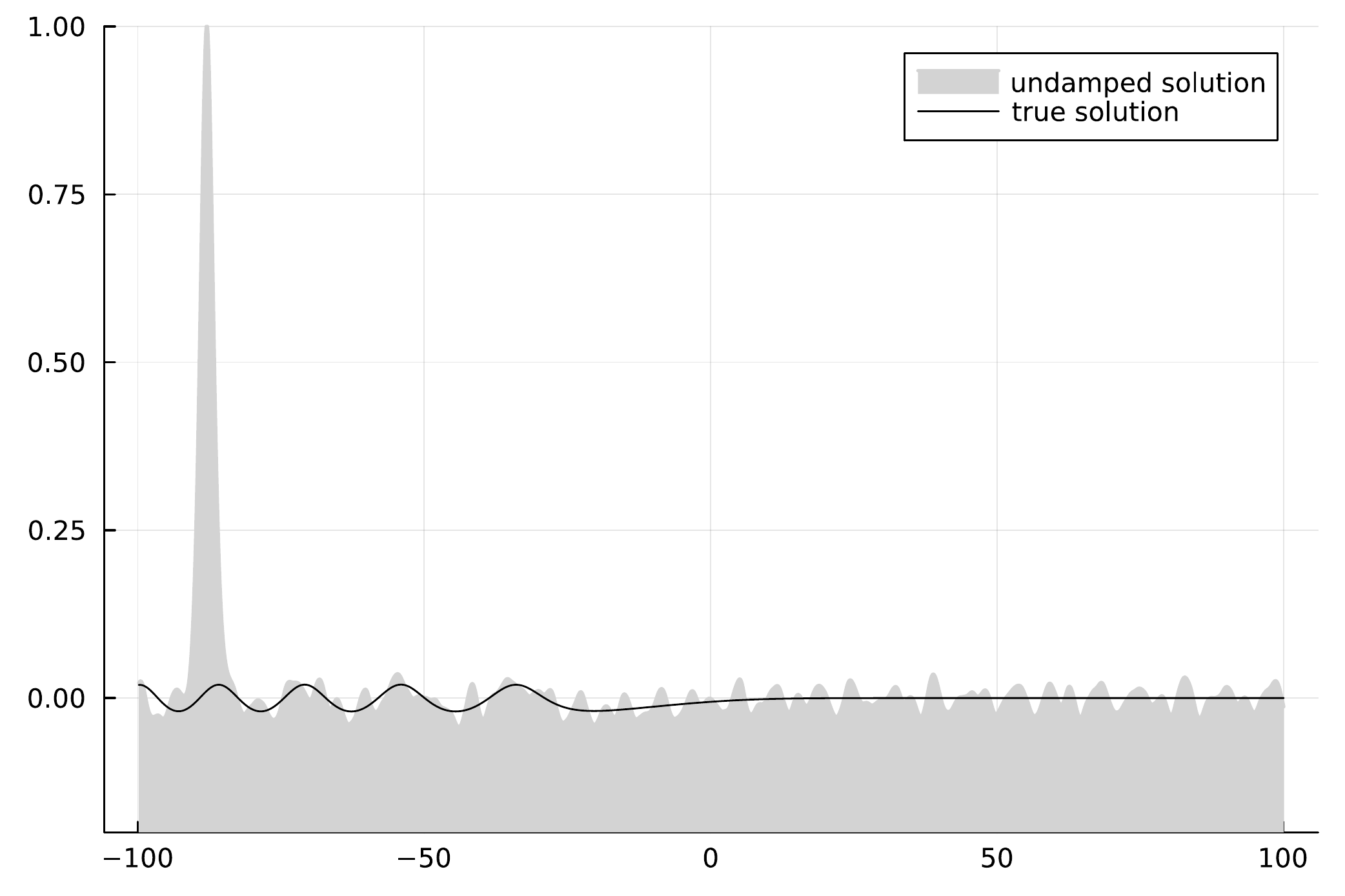}
  \caption{Undamped solution to the KdV equation computed on \([-200,200]\) with \(2^{10}\) grid points.}
  \label{fig:fig3}
\end{figure}

\begin{figure}[h!]
    \centering
    \includegraphics[width = 0.8\linewidth]{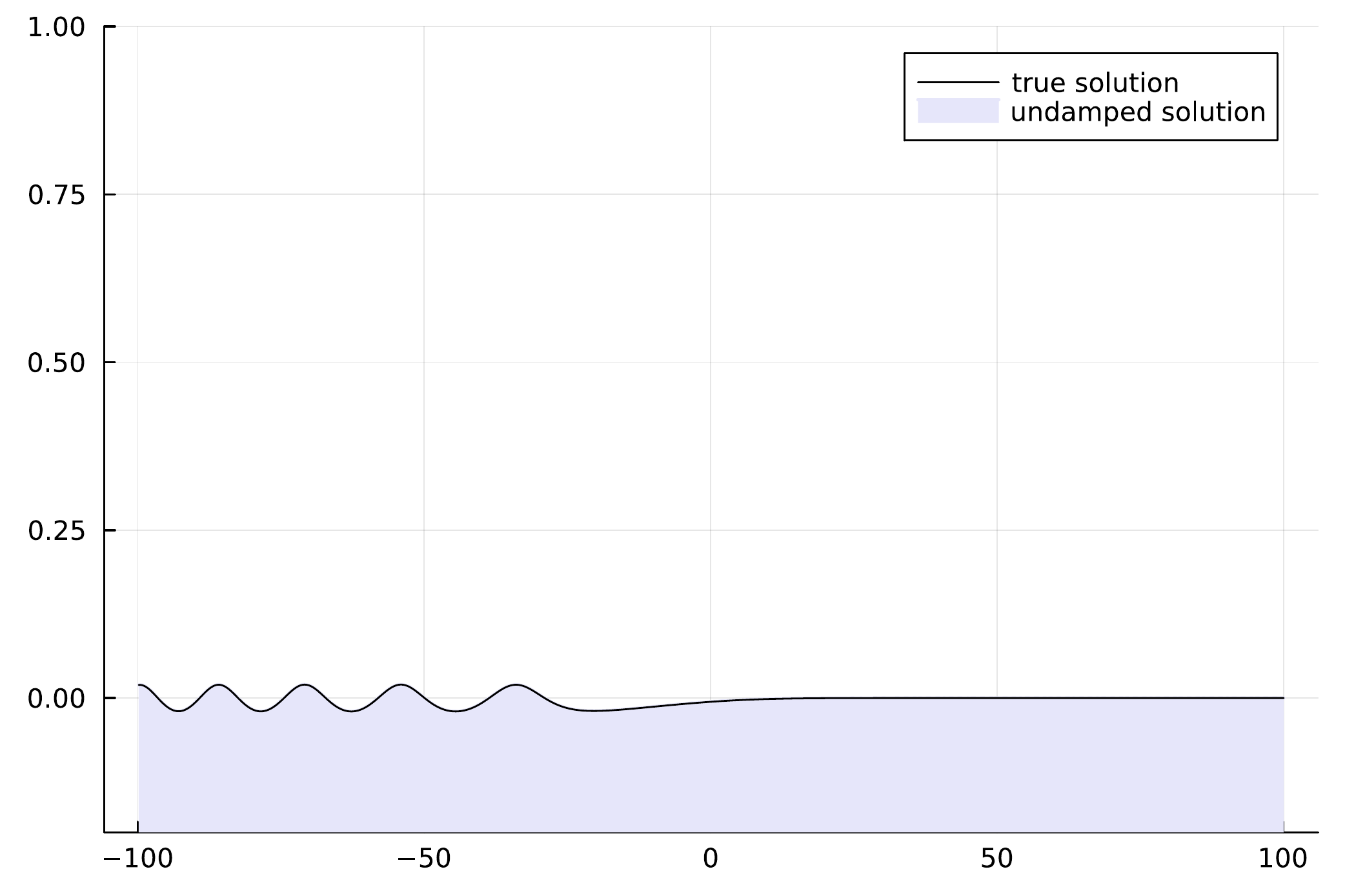}
    \caption{Undamped solution to the KdV equation computed on \([-2500,2500]\) with \(2^{14}\) grid points.}
    \label{fig:udgood}
\end{figure}

Without artificial damping, the primary way to mitigate this problem is to increase the size of the interval of computation. This improves the accuracy of the solution by creating a buffer between the periodic instances of the solution. However, this means that observing the long-term behavior of a solution requires a sizable interval and thus a large runtime. In this section, we compute an accurate solution to \eqref{eq:kdv} with \(q(x,0) = 1.3\ee^{-x^2}\) at \(t = 150\).\footnote{All solutions in this section are computed with a time step of 0.01.} We compare it to the true solution\footnote{The solution used as the ground truth in this section comes from \cite{trogdonkdv}.} on \([-99.85,100.05]\), seeking a maximum error on the order of \(10^{-8}\). To demonstrate the effect of the periodic Fourier method on the accuracy of the solution, we use an interval of computation, \([-200,200]\), as shown in Figure~\ref{fig:fig3}.
This computation takes 13 seconds and produces a maximum error of 0.99. We can clearly see that the Fourier method has introduced oscillations and a soliton that do not appear in the true solution. In order to achieve a more accurate solution, we increase the interval of computation to \([-10000,10000]\), as shown in Figure~\ref{fig:udgood}.
This solution still falls short of a reasonable accuracy goal of of \(4\cdot 10^{-6}\). Computing it took 869 seconds (\(\approx 14.5\) minutes). Computing a solution on a larger interval, at a later time, or with a better error can become prohibitively time-consuming to do on most computers. For this reason, we introduce artificial damping.

\subsection{Artificially-damped solutions of the KdV equation}
\label{sec:damping}

As in Section~\ref{sec:damping1}, we solve a modified version of the KdV equation \eqref{eq:kdv},
\begin{equation*}
    q_t+6qq_x+q_{xxx}= k_1\left(\sigma(x)q_{x}\right)_x-k_2(1-\gamma(x))q.
\end{equation*}
We damp the small-amplitude dispersive tail on the left side of the interval by solving the heat equation with Strang-splitting. We choose \(k_1 = 1\) and 
\begin{equation}
    \sigma(x) = \sigma(x;\ell_1,\ell_2) = 1 - \frac{1}{2}\left(\tanh(x-\ell_1)+1\right) + \frac{1}{2}\left(\tanh(-x-\ell_2)+1)\right),
    \label{eq:sigma}
\end{equation}
where \(\ell_1 = -L+\frac{L}{2}-10\) and \(\ell_2 = L-5\). We perform this diffusion at every time step (\(f_1 = 1\)).
\begin{figure}
    \begin{minipage}{0.5\textwidth}
        \centering
        \includegraphics[width=0.9\textwidth]{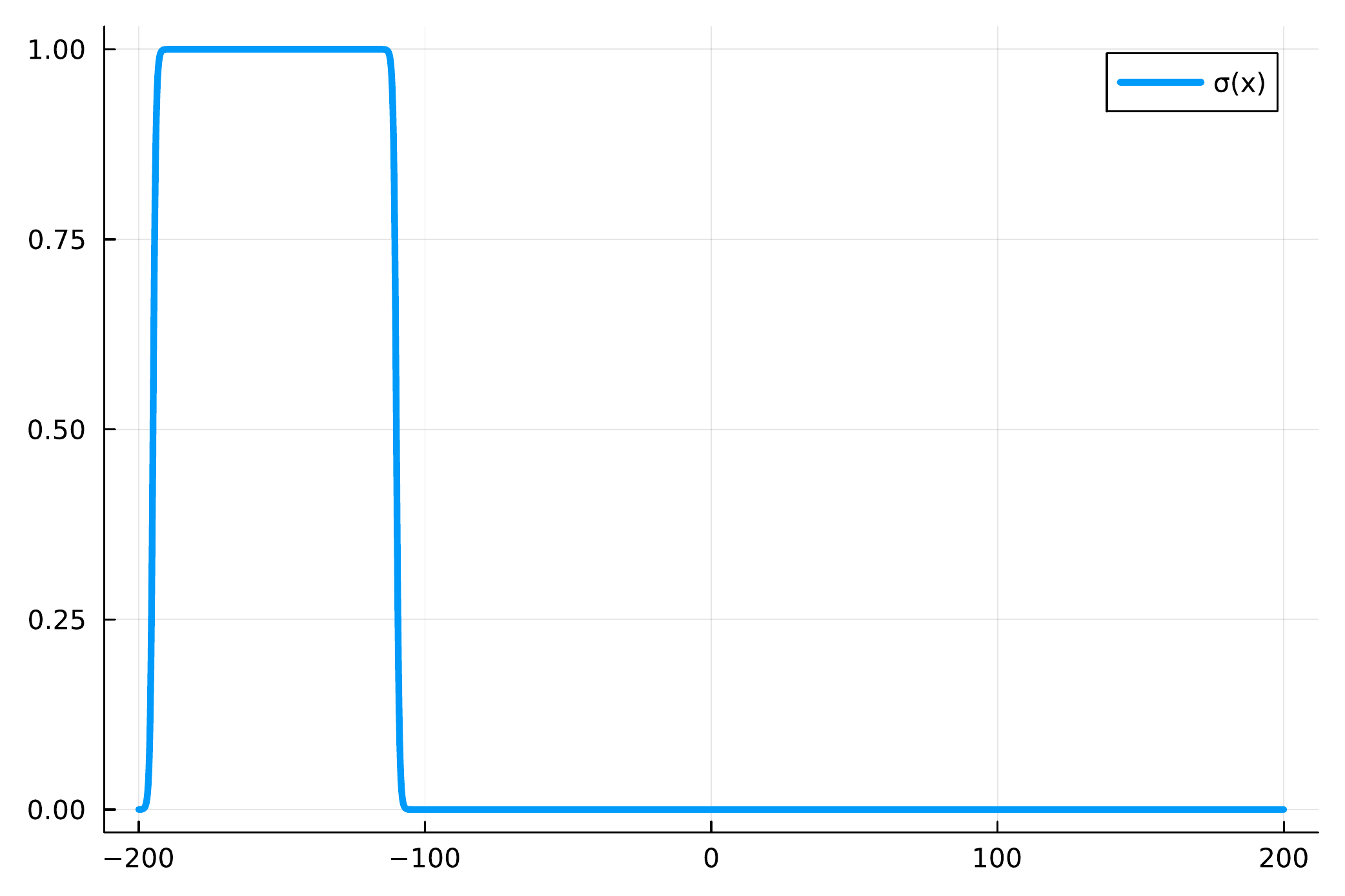} % first figure itself
        \caption*{\footnotesize{Heat eq. damping coefficient, \(\sigma(x)\).}}
    \end{minipage}\hfill
    \begin{minipage}{0.5\textwidth}
        \centering
        \includegraphics[width=0.9\textwidth]{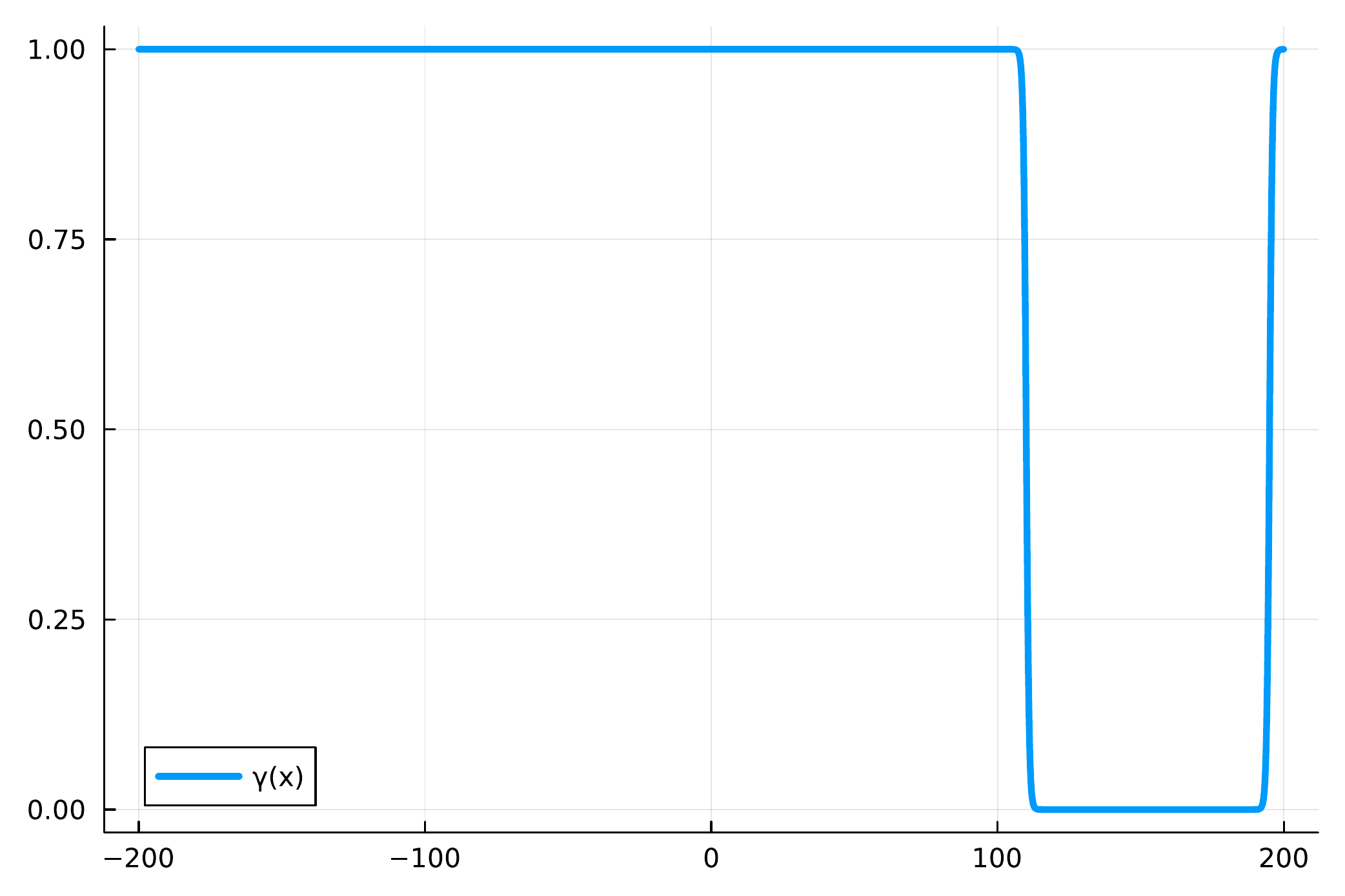} % second figure itself
        \caption*{\footnotesize{Exponential decay damping factor, \(\gamma(x)\).}}
    \end{minipage}
\end{figure}
The soliton on the right side of the interval requires a more aggressive damping technique. We simulate rapid exponential decay as \(k\rightarrow\infty\) by multiplying the solution values by
\begin{equation}
    \gamma(x)=1-\sigma(-x)
    \label{eq:gamma}
\end{equation}
every 1000 time steps \((f_2 = 1000)\). We compute an artificially-damped solution on \([-600,600]\) with \(2^{12}\) grid points and a time step of 0.01. This solution produces errors that are better than the undamped solution in Figure~\ref{fig:udgood}, its error is \(4\cdot 10^{-8}\); however, it only takes 157 seconds (\(\approx 2.6\) minutes) to compute. To summarize, we have achieved an error that is 100 times smaller in \(18\%\) of the runtime. Further comparisons of solutions computed with and without damping are in Table~\ref{tab:tab1}.
\begin{table}
      \centering
      \begin{tabular}{l r r c}
        \toprule
        \textbf{Parameters} & \textbf{Undamped Error} & \textbf{Damped Error}\\
        \midrule
        \(L = 100,\;m = 2^9\) & \(1.05\) (8 sec) & \(0.03\) (20 sec) \\
        \(L = 200,\;m = 2^{10}\) & 0.99 (13 sec) & 0.0006 (65 sec)\\
        \(L = 600,\; m = 2^{12}\) & 0.01 (52 sec) & \(4\cdot 10^{-8}\) (157 sec)\\
        \(L = 1200,\; m = 2^{13}\) & 0.002 (107 sec) & \(4\cdot 10^{-8}\) (284 sec)\\
        \(L = 2500,\; m = 2^{14}\) & 0.0002 (140 sec) &  \(\cdot\)\\
        \(L = 5000,\; m = 2^{15}\) & \(0.0001\) (305 sec) &  \(\cdot\)\\
        \(L  = 10000, \; m = 2^{16}\) & \(4\cdot 10^{-6}\) (869 sec) & \(\cdot\)\\
        \(L = 20000, \; m = 2^{17}\) & \(1\cdot 10^{-6}\) (2812 sec) & \(\cdot\)\\
        \(L = 40000, \; m = 2^{18}\) & \(1\cdot 10^{-6}\) (4260 sec) & \(\cdot\)\\
        \bottomrule
      \end{tabular}
      \caption{Solutions evaluated on an interval \([-L,L]\) with \(m\) grid points.}
      \label{tab:tab1}
    \end{table}
We see that the errors for the damped solutions saturate on the order of \(10^{-8}\) as the size of the interval increases. This is because we are using a fourth-order time-stepping method (RK4) and a time step of 0.01. 
    
\subsection{A note on damping parameters}

\begin{figure}[h!]
    \centering
    \includegraphics[width = \linewidth]{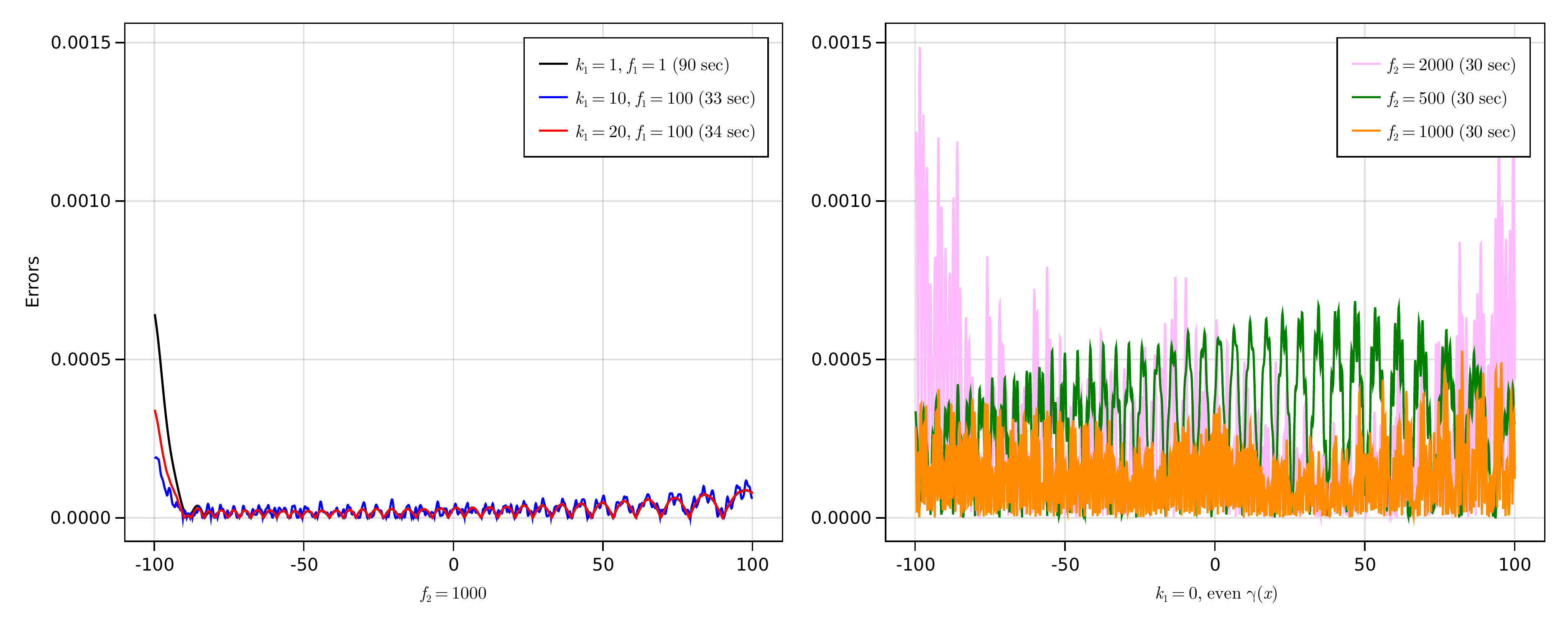}
    \caption{Errors for solutions to \eqref{eq:kdv} at \(t = 150\), computed on \([-200,200]\) with \(m = 2^{10}\) grid points and a time step of 0.01. Error is defined as the absolute difference between the true solution and the artificially-damped solution at each grid point.}
    \label{fig:diffparams}
\end{figure}

The choices of \(k_1\), \(k_2\), \(f_1\), \(f_2\), \(\sigma(x)\), and \(\gamma(x)\) will determine the severity of damping and its region. We always consider \(k_2\rightarrow\infty\) to simulate rapid exponential decay, but the other parameters can be chosen judiciously. For this problem, we found that setting \(k_1 = 1\), \(f_1 = 1\), and \(f_2 = 1000\), as well as choosing \(\sigma(x)\) to be non-zero on the left and \(\gamma(x)\) to be zero on the right, produced an accurate solution on a smaller interval. However, it is important to note that different choices are valid and may be advantageous depending on the problem. When using our method, experimentation with these parameters is encouraged to determine values that best suit the specific problem. In Figure~\ref{fig:diffparams}, we demonstrate how the error is affected by the damping parameters.
The first plot shows the effect of altering \(k_1\) and \(f_1\), while holding \(f_2\) constant and implementing \(\gamma(x)\) damping every 1000 time steps. In the second plot, we do not use \(\sigma(x)\) damping and instead use an even \(\gamma(x)\), as shown in Figure~\ref{fig:fig8}. We alter \(f_2\) to demonstrate the importance of choosing an optimal value. Several observations are in order. 

\begin{enumerate}
    \item There exists an important trade-off between accuracy and runtime of \(\sigma(x)\) and \(\gamma(x)\) damping, where use of only \(\gamma(x)\) damping is faster, but including \(\sigma(x)\) damping improves errors.
    \item There appears to be a balance between frequency and amplitude of \(\sigma(x)\) and \(\gamma(x)\) damping that maintains a shorter runtime while achieving high accuracy. 
    \item Different choices of parameters result in different oscillatory frequencies of the errors across the computational domain. It may be advantageous to choose parameters that display less frequent oscillations.
    \item In the left plot, errors spike at the left side of the computational domain. This could be due to the high-frequency dispersive tail that is difficult to approximate numerically, or the rightward traveling soliton that has wrapped around to the left side of the domain. Note that the blue curve corresponding to the most aggressive \(\sigma(x)\) damping on the plot minimizes this spike in the error. 
\end{enumerate}

\subsection{On the linearized Korteweg-de Vries equation}

\begin{figure*}
    \centering
    \includegraphics[scale = 0.5]{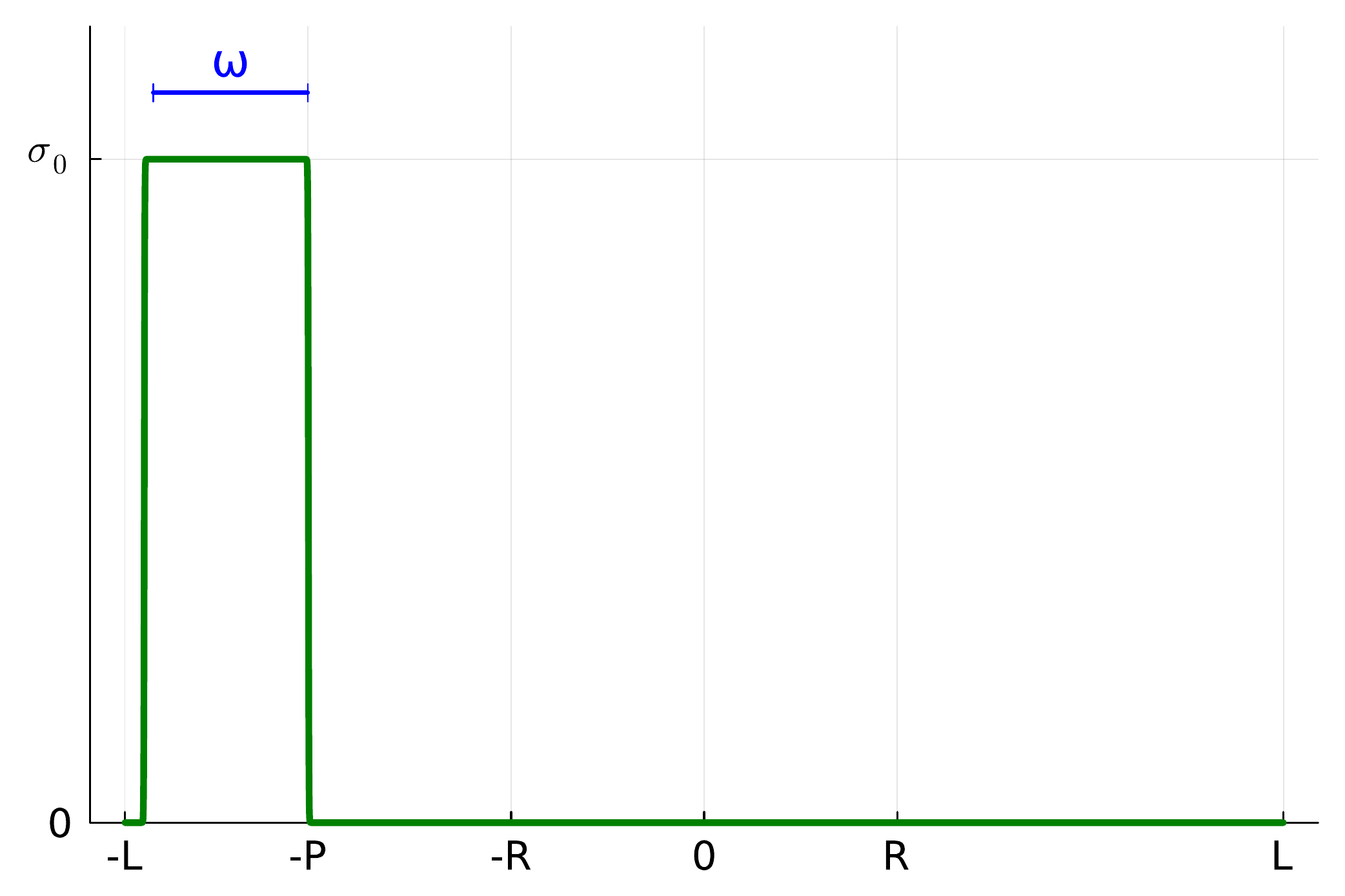}
    \caption{General diffusion coefficent \(\sigma(x)\).}
    \label{fig:linearizedsigma}
\end{figure*}

We consider the linearized Korteweg-de Vries equation,
\begin{equation}
    q_t+q_{xxx}=0,
    \label{eq:linkdv}
\end{equation}
in order to more concretely, but still heuristically, demonstrate the soundness of our numerical damping method. We introduce damping by considering a modified equation,
\begin{equation*}
    \Tilde{q}_t+\Tilde{q}_{xxx}=(\sigma(x)\Tilde{q}_{x})_x,
\end{equation*}
where again \(\sigma(x)\) is a non-zero diffusion coefficient that determines the region of damping, as shown in \eqref{fig:linearizedsigma}.

We consider an initial condition that is a wave packet with fundamental frequency \(k\) centered at the origin. This wave packet travels with approximate velocity \(-3k^2\) and at time \(T = \frac{P}{k^2}\), it arrives at\footnote{Here, we change the notation from \(P_-\) to \(P\) for simplicity and because we are only performing damping on the left, so distinction between \(P_-\) and \(P_+\) is no longer necessary.} \(x=-P\) and damping begins. It remains in the damping region for a time \(\tau=\frac{\omega}{k^2}\) and experiences damping of magnitude \(\ee^{-\sigma_0k^2\tau}=\ee^{-\sigma_0\omega}=:\varepsilon_0\). We know that the true solution to the linearized KdV equation can be expressed as
\begin{equation}
    q(x,t)=\frac{1}{2\pi}\int_{-\infty}^\infty \ee^{\ii kx+\ii k^3t}\hat{q}_0(k)\dd k
    \label{eq:linkdvsol}
\end{equation}
where \(\hat{q}_0(k)\) is the Fourier transform initial condition:
\[\hat{q}_0(k)=\int_{-\infty}^{\infty}\ee^{-\ii kx}q(x,0)\dd x.\] 
While many calculations below are general, the final conclusions will be valid for the choice
\begin{align*}
\hat q_0(k) = \ee^{-k^2}.
\end{align*}

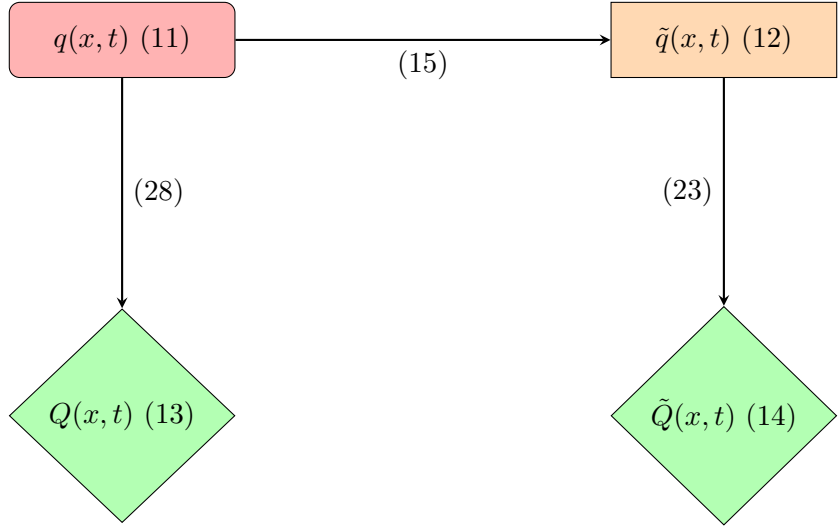
\begin{figure}[htp]
\centering
\begin{tikzpicture}[node distance=4cm]
\node (start) [startstop] {\(q(x,t)\) \eqref{eq:linkdvsol}};
\node (pro1) [process, right of=start, xshift = 4cm] {\(\Tilde{q}(x,t)\) \eqref{eq:dampintegral}};
\node (dec1) [decision, below of=pro1, yshift = -1cm] {\(\Tilde{Q}(x,t)\) \eqref{eq:heuristic}};
\node (dec2) [decision, below of=start, yshift = -1cm] {\(Q(x,t)\) \eqref{eq:capitalQ}};
\draw [arrow] (start) -- node[anchor=north] {\eqref{eq:wholedampedapprox}} (pro1);
\draw [arrow] (start) -- node[anchor=west] {\eqref{eq:capitalQbound}} (dec2);
\draw [arrow] (pro1) -- node[anchor=east] {\eqref{eq:mostimportantbound}}(dec1);
\end{tikzpicture}
\caption{In this section, we compare the undamped \eqref{eq:capitalQ} and damped \eqref{eq:heuristic} periodic approximations of the true solution. The damped approximation is a periodic approximation of the damped solution on the whole line \eqref{eq:dampintegral}.}
\label{fig:flowchart}
\end{figure}
We are led to conjecture that the damped solution can be approximated as 
\begin{equation}
    \Tilde{q}(x,t)\approx \frac{1}{2\pi}\int\displaylimits_{k:\;|3k^2t|<P}\ee^{\ii kx+\ii k^3t}\hat{q}_0(k)\dd k+\frac{\varepsilon_0}{2\pi}\int\displaylimits_{k:\;|3k^2t|>P}\ee^{\ii kx+\ii k^3t}\hat{q}_0(k)\dd k.
    \label{eq:dampintegral}
\end{equation}
Consider the periodic problem,
\begin{align*}
    Q_t+Q_{xxx}&=0, \quad Q(x + 2L,t) = Q(x,t),\quad
    Q(x,0) \approx q_0(x), \quad x \in [-L,L].
\end{align*}
We construct a solution by discretizing the integral in \eqref{eq:linkdvsol}:
\begin{equation}
    Q(x,t)=\frac{1}{2L}\sum_{j=-\infty}^\infty \ee^{\ii kx+\ii k^3t}\hat{q}_0(k_j), \quad k_j=\frac{\pi}{L}j.
    \label{eq:capitalQ}
\end{equation}
From this, we might expect that we can approximate a damped solution of
\begin{align*}
    \Tilde{Q}_t+\Tilde{Q}_{xxx}&=\sigma(x)\Tilde{Q}_{xx}\quad \tilde Q(x + 2L,t) = \tilde Q(x,t),\\
    \tilde Q(x,0) &\approx q_0(x), \quad x \in [-L,L].
\end{align*}
using 
\begin{equation}
    \Tilde{Q}(x,t)\approx \frac{1}{2L}\sum_{\substack{j=-\infty \\ |3k_j^2t|<P}}^\infty \ee^{\ii k_jx+\ii k_j^3t}\hat{q}_0(k_j) + \frac{\varepsilon_0}{2L}\sum_{\substack{j=-\infty \\ |3k_j^2t|>P}}^\infty \ee^{\ii k_jx+\ii k_j^3t}\hat{q}_0(k_j).
    \label{eq:heuristic}
\end{equation}

This provides a heuristic for us to analyze what the damping method is doing. However, we must first confirm that \eqref{eq:heuristic} is indeed a good approximation of what our numerical method produces. We compute solutions to the linearized KdV equation at \(T=150\) using \eqref{eq:heuristic} and our numerical method and measure their maximum difference, as shown in \eqref{tab:tab2}. We choose \(\hat{q}_0(k)=\ee^{-k^2}\) and define \(\sigma(x)\) as in \eqref{eq:sigma}. We use a time step of 0.01 for the numerical solutions. 

\begin{table}[H]
      \centering
      \begin{tabular}{l r r c}
        \toprule
        \textbf{Parameters} & \textbf{Error} \\
        \midrule
        \(L = 100,\;m = 2^9\) &  \(0.02\)\\
        \(L = 200,\;m = 2^{10}\) & \(0.01\)\\
        \(L = 600,\; m = 2^{12}\) & \(0.007\) \\
        \(L = 1200,\; m = 2^{13}\) & \(0.003\)\\
        \bottomrule
      \end{tabular}
      \caption{Solutions are evaluated on an interval \([-L,L]\) with \(m\) grid points. Error is calculated as the maximum absolute difference between the numerical method and \eqref{eq:heuristic} on the computational domain, excluding the damping region.}
      \label{tab:tab2}
    \end{table}

From the table, we argue that \eqref{eq:heuristic} is a reasonable heuristic for what our numerical damping method is producing. We begin analysis of how well \eqref{eq:heuristic} approximates the true solution \eqref{eq:linkdvsol} by first considering the damped approximation on the whole line \eqref{eq:dampintegral} as an approximation of the true solution \eqref{eq:linkdvsol}:
\begin{align*}
    &|\eqref{eq:linkdvsol}-\eqref{eq:dampintegral}| \\
    &= 
    \left|\frac{1}{2\pi}\int_{-\infty}^\infty \ee^{\ii kx+\ii k^3t}\hat{q}_0(k)\dd k - \frac{1}{2\pi}\int\displaylimits_{k:\;|3k^2t|<P}\ee^{\ii kx+\ii k^3t}\hat{q}_0(k)\dd k-\frac{\varepsilon_0}{2\pi}\int\displaylimits_{k:\;|3k^2t|>P}\ee^{\ii kx+\ii k^3t}\hat{q}_0(k)\dd k\right|
    \\& \leq \frac{(1-\varepsilon_0)}{2\pi}\int\displaylimits_{k:\;|3k^2t|>P}\left|\ee^{\ii kx+\ii k^3t}\hat{q}_0(k)\right|\dd k \quad \leq \quad\frac{(1-\varepsilon_0)}{2\pi}\int\displaylimits_{k:\;|3k^2t|>P}\left|\hat{q}_0(k)\right|\dd k.
\end{align*}
Since \(\varepsilon_0\) is small, \((1-\varepsilon_0)\) will not be. Therefore, we need the integral in
\begin{equation}
    \frac{(1-\varepsilon_0)}{2\pi}\int\displaylimits_{k:\;|3k^2t|>P}\left|\hat{q}_0(k)\right|\dd k
    \label{eq:wholedampedapprox}
\end{equation}
to be small to have a reasonable approximation. Throughout the rest of this section, we assume that \(|\hat{q_0}(k)|\leq \ee^{-k^2}\), and use the rough estimate
\[
\int\displaylimits_{k:\;|3k^2t|>P}\ee^{-k^2}\dd k \leq \int_{\sqrt{\frac{P}{3t}}}^{\infty}\ee^{-k}\dd k.
\]
If we desire that this contribution is bounded by \(\varepsilon_1\), we arrive at the following condition:
\begin{equation}
    \sqrt{\frac{P}{3t}} \geq \ln\left(\frac{1}{\varepsilon_1}\right)
    \label{eq:PTcondition}
\end{equation}
Next, we consider the difference between the damped periodic approximation \eqref{eq:heuristic} and the damped whole-line approximation \eqref{eq:dampintegral}. We begin by considering whether 
\begin{equation}
    \frac{1}{2L}\sum_{\substack{j=-\infty \\ |3k_j^2t|<P}}^\infty \ee^{\ii k_jx+\ii k_j^3t}\hat{q}_0(k_j)
    \label{eq:riemannsum}
\end{equation}
is a good approximation of 
\begin{equation*}
    \frac{1}{2\pi}\int\displaylimits_{k:\;|3k^2t|<P}\ee^{\ii kx+\ii k^3t}\hat{q}_0(k)dk
\end{equation*}
(i.e. whether the first term in \eqref{eq:heuristic} is a good approximation of the first term in \eqref{eq:dampintegral}). Since the second terms are both multiplied by \(\varepsilon_0\), which is small, we ignore the effect of these terms in the comparison. We are seeking a bound for 
\begin{equation}
\left|\frac{1}{2\pi}\int_{-\sqrt{\frac{P}{3t}}}^{\sqrt{\frac{P}{3t}}}\ee^{\ii kx+\ii k^3t}\hat{q}_0(k)\dd k-\frac{1}{2L}\sum_{j=-M}^M \ee^{\ii k_jx+\ii k_j^3t}\hat{q}_0(k_j)\right|.
\label{eq:errorbound2}
\end{equation}

We first note that the grid of \(k_j\) values may not line up perfectly such that \(k_{-M}=-\sqrt{\frac{P}{3t}}\) and \(k_M=\sqrt{\frac{P}{3t}}\). Generally, \(k_{-M}\) will be the largest negative \(k_j\) value that is greater than \(-\sqrt{\frac{P}{3t}}\) and \(k_M\) will be the largest \(k_j\) that is less than \(\sqrt{\frac{P}{3t}}\). Then, the integral that the Riemann sum is actually approximating is
\begin{equation}
\frac{1}{2\pi}\int_{k_{-M}}^{k_M}\ee^{\ii kx+\ii k^3t}\hat{q}_0(k)\dd k.
\label{eq:actualbounds}
\end{equation}
This differs from the integral we want to approximate by at most 
\begin{equation}
    \frac{2\pi}{L}\max\limits_{|k|>\sqrt{\frac{P}{3t}}-\frac{\pi}{L}}|\hat{q}_0(k)|,
    \label{eq:differentbounds}
\end{equation}
since \(k_{-M}\) and \(k_M\) cannot differ from the integral bounds by more than the size of the grid spacing, \(\frac{\pi}{L}\). This is one component of the error in \eqref{eq:errorbound2} and, because of the condition on \(P\) and \(t\) \eqref{eq:PTcondition}, we assume this contribution is small and ignore its effect on the overall error bound. The second component of the error is more significant and comes from how well the Riemann sum \eqref{eq:riemannsum} approximates \eqref{eq:actualbounds}. \\

\begin{theorem}[The exponentially convergent trapezoidal rule for periodic integrals \cite{expconvergent}]\label{thm:T_thm}
  An integral,
\begin{equation}
    I = \int_{0}^T v(\theta)\dd\theta,
    \label{eq:specialintegral}
\end{equation}
can be approximated by a trapezoidal rule approximation,
\begin{equation*}
    I_N = \frac{2\pi}{N}\sum_{k=1}^Nv(\theta_k),
\end{equation*}
with an error bound,
\begin{equation}
    \left|I - I_N\right|\leq \frac{2TA}{\ee^{\frac{2\pi A N}{T}}-1},
    \label{eq:fancybound}
\end{equation}
for any \(N\geq 1\) if \(v(\theta)\) is T-periodic, analytic, and satisfies \(|v(\theta)|\leq A\) in the strip \(-a< \imag \theta<a\) for some \(a>0\). 
\end{theorem}

We use \eqref{thm:T_thm} on the integral in \eqref{eq:actualbounds} by assuming that \(\hat{q}_0(k)\rightarrow 0\) as \(|k|\rightarrow \infty\) sufficiently quickly such that it is approximately periodic on \([k_{-M},k_M]\). Defining \(\theta = k+k_M\) and assuming \(k_{-M}\approx -k_M\) puts \eqref{eq:actualbounds} in the desired form of \eqref{eq:specialintegral} with
\begin{equation*}
    v(\theta) = \frac{1}{2\pi}\ee^{\ii(\theta-k_M)x+\ii(\theta-k_M)^3t}\hat{q_0}(\theta-k_M)
\end{equation*}
and \(T = 2k_M\). Since we have assumed \(|\hat{q_0}(k)| \leq \ee^{-k^2}\), this allows us to compute a specific upper bound for \(|v(\theta)|=|v(E+\ii\eta)|\) for \(|\eta|<a\):
\begin{equation*}
    |v(\theta)|\leq \frac{1}{2\pi}\ee^{aR+a^3t+a^2}\ee^{(3\eta t-1)(E-k_M)^2}
\end{equation*}
This leads us to impose that \(a = \frac{1}{3t}\) such that \(3\eta t - 1 < 0\) and the exponential is bounded. From this, we conclude that
\[
\ee^{(3\eta t-1)(E-k_M)^2} \leq 1
\]
and arrive at an upper bound on \(|v(\theta)|\):
\begin{equation*}
    A = \frac{1}{2\pi}\ee^{aR+a^3 t+a^2},
\end{equation*}
where \(a = \frac{1}{3t}\). Now we apply \eqref{eq:fancybound}, with \(T = 2k_M\) and \(N = \frac{2k_ML}{\pi}\), to arrive at the following error bound:
\begin{equation}
    |\eqref{eq:actualbounds}-\eqref{eq:riemannsum}|\lesssim \frac{2k_M}{\pi}\ee^{aR+a^3t+a^2-2aL}.
    \label{eq:mostimportantbound}
\end{equation}
Assuming \(P\) is sufficiently large according to \eqref{eq:PTcondition} and \(\varepsilon_0\) is small, we ignore the other smaller contributions to the error from \eqref{eq:differentbounds} and the terms multiplied by \(\varepsilon_0\) in \eqref{eq:linkdvsol} and \eqref{eq:heuristic}. We consider \eqref{eq:mostimportantbound} to be an upper bound on \(|\eqref{eq:linkdvsol}-\eqref{eq:heuristic}|\). 
If an error of \(\epsilon\) is desired, then \(L\) and \(t\) must satisfy
\begin{equation}
    L\geq \frac{3}{2}t\ln\left(\frac{1}{\epsilon}\right)+\frac{3}{2}t\ln\left(\frac{2}{\pi}\sqrt{\frac{P}{3t}}\right)+\frac{1}{2}R+\frac{2}{9t},
    \label{eq:periodicinterval}
\end{equation}
where we have used the fact that \(k_M\approx \sqrt{\frac{P}{3t}}\). 

Our goal is to make a theoretical comparison between damped and undamped solutions, so we consider
\begin{equation}
    |\eqref{eq:linkdvsol}-\eqref{eq:capitalQ}|=
    \left|\frac{1}{2\pi}\int_{-\infty}^\infty \ee^{\ii kx+\ii k^3t}\hat{q_0}(k)\dd k - \frac{1}{2L}\sum_{j=-\infty}^{\infty}\ee^{\ii k_j x + \ii k_j^3 t}\hat{q_0}(k_j)\right|,
    \label{eq:undampedbound}
\end{equation}
an upper bound on the error for the undamped periodic approximation of the solution. 

\begin{theorem}[The exponentially convergent trapezoidal rule for whole-line integrals \cite{expconvergent}]\label{thm:bigthm}
An integral,
\begin{equation*}
    I = \int_{-\infty}^\infty \omega(k)\dd k,
\end{equation*}
can be approximated by a trapezoidal rule approximation,
\begin{equation*}
    I_h^{[n]} = h\sum_{j=-M}^M\omega(k_j),
\end{equation*}
with an error bound,
\begin{equation}
    \left|I - I_h^{[n]}\right|\leq \frac{2A}{\ee^{\frac{2\pi a}{h}}-1},
    \label{eq:fancybound_line}
\end{equation}
if \(\omega(k)\) is analytic in the strip \(-a < \mathrm{Im}\, k <a\) for some \(a>0\), \(\omega(k)\rightarrow 0\) uniformly as \(|k|\rightarrow \infty\) in the strip, and for some \(A\),
\begin{equation*}
    \int_{-\infty}^{\infty}|\omega(E+\ii\eta)|\dd E\leq A \quad \text{for all} \quad \eta \in (-a,a).
\end{equation*}
\end{theorem}

We use \eqref{thm:bigthm} by defining \(\omega(k) = \frac{1}{2\pi}\ee^{\ii kx+\ii k^3t}\hat{q_0}(k)\) and \(h = \triangle k = \frac{\pi}{L}\) to derive an upper bound for \eqref{eq:undampedbound}.
We begin by seeking an upper bound on the integral of \(\omega(E+\ii\eta)\):
\begin{align*}
    \frac{1}{2\pi}\int_{-\infty}^\infty \left|\ee^{\ii(E+\ii\eta)x}\ee^{\ii(E+\ii\eta)^3t}\hat{q_0}(E+\ii\eta)\right|\dd E & \leq \frac{1}{2\pi}\int_{-\infty}^\infty \left|\ee^{\eta|x|}\ee^{(\eta^3+3E^2\eta)t}\hat{q_0}(E+\ii\eta)\right| \dd E\\ & \leq \frac{1}{2\pi}\ee^{aR+a^3t}\int_{-\infty}^{\infty}\left|\ee^{3\eta tE^2}\hat{q_0}(E+\ii\eta)\right|\dd E.
\end{align*}
With \(|\hat{q}_0(k)|\leq \ee^{-k^2}\), we have that 
\[
A = \frac{C}{2\pi}\ee^{aR+a^3t+a^2}
\]
where 
\begin{equation}
    C = \int_{-\infty}^{\infty}\ee^{(3\eta t - 1)E^2}\dd E = \sqrt{\frac{\pi}{1-3\eta t}} \quad \text{for} \quad \eta < \frac{1}{3t}.
    \label{eq:C_bound}
\end{equation}
We impose that \(|\eta| < a = \frac{1-\delta}{3t}\) for $0 < \delta < 1$ to ensure that 
\(C\) is finite. Now, we obtain an error bound of the form \eqref{eq:fancybound_line},
\begin{equation}
    |\eqref{eq:linkdvsol}-\eqref{eq:capitalQ}| \lesssim \frac{C}{\pi}\ee^{aR+a^3t+a^2-2aL}
    \label{eq:capitalQbound}
\end{equation}
which we use to derive a similar relationship between \(L\) and \(t\) for a desired error of \(\epsilon\),
\begin{equation}
    L\geq \frac{3}{2(1-\delta)}t\ln\left(\frac{1}{\epsilon}\right)+\frac{3}{2(1-\delta)}t\ln\left(\frac{1}{\pi}C\right)+\frac{1}{2(1-\delta)}R+\frac{2}{9t(1-\delta)}.
    \label{eq:wholeline}
\end{equation}

Note that this differs from the damped error bound \eqref{eq:periodicinterval} in that all of the coefficients are larger by a factor of \(\frac{1}{1-\delta}\) and that in the second term we previously had
\(
\ln\left(\frac{2}{\pi}\sqrt{\frac{P}{3t}}\right)
\)
whereas now we have
\(
\ln\left(\frac{1}{\pi}C\right)
\). The value of \(C\) can be decreased by increasing \(\delta\), however this also increases the coefficients of every term in the bound. On the other hand, \(P\) is a damping parameter that, although it must be sufficiently large according to \eqref{eq:PTcondition}, can be chosen somewhat judiciously. Through our analysis of the damped and undamped errors, the benefit of artificial damping is made clearer. The introduction of the additional parameter \(P\) allows us to impose a less restrictive relationship on the size of the computational domain, \(L\), as it relates to the final time, \(t\). This provides a partial explanation for our observations when solving the nonlinear Korteweg-de Vries equation in \eqref{eq:kdv}, which is that artificially-damped solutions achieved a given error with a smaller computational domain than undamped solutions. However, we believe that this theoretical analysis is pessimistic in its projection of the advantages of our damping method. The numerical data suggests far less restrictive conditions on \(L\) are required by the artificial damping method than the small improvement predicted by the error bound. 

One possible indication of this analysis, something that might hint at a ceiling on the improvements of the proposed method, is that the lower bounded for $L$ in both cases are, to leading order, linear in $t$.  So, it is likely that introducing damping does not change this order, keeping it linear, but introducing damping reduces the leading constant.

\begin{comment}
\textbf{A lower bound for the undamped solution }
Errors propagate in undamped solutions of the KdV equation when the dispersive tail flows out of the computational domain from the left and reenters it on the right. In an ideal situation, we would have the true solution up to a time \(t-\triangle t\) and would only need the undamped Fourier method to move forward one time step to obtain the solution at the desired time \(t\). In this case, the error would be at least as bad as the magnitude of the solution on the leftside of the compuational domain, \(q(-L,t)\). From (*insert citation here*), an asymptotic formula for the solution of the linearized KdV equation is given by
\begin{equation}
    q(x,t)=\frac{|\hat{q_0}(z_0)|}{\sqrt{3\pi t z_0}}\cos\left(2tz_0^3-\frac{\pi}{4}-\textrm{arg}(f(z_0))\right)+\mathcal{O}(t^{-\frac{3}{2}}),
\end{equation}
where \(z_0=\sqrt{\frac{|x|}{3t}}\).
Therefore, requiring that this lower bound be less than \(\varepsilon\),
\begin{equation*}
    |q(-L,t)| \approx \frac{|\hat{q_0}(z_0)|}{\sqrt{3\pi t z_0}} = C_1\frac{|\hat{q_0}(z_0)|}{(tL)^{\frac{1}{4}}} < \varepsilon,
\end{equation*}
gives a condition on the relationship between \(L\) and \(t\):
\begin{equation}
    L\geq 3t\ln\left(\frac{1}{\varepsilon}\right)-\frac{3}{4}t\ln(tL)-3t\ln(C_1).
    \label{eq:undampedLcondition}
\end{equation}
\end{comment}

\section{Nonlinear Schrödinger equation}
\label{sec:NLS}

In this section, we compute solutions to the nonlinear Schrödinger (NLS) equation \cite{ognls},
\begin{equation}
    \ii q_t+q_{xx}+2|q|^2q=0,
    \label{eq:nls}
\end{equation}
with an initial condition, \(q(x,0)=(1+x)\ee^{\ii x-0.7x^2}\). The real part of the solution has small-amplitude oscillations traveling in both directions.

\subsection{The Fourier method for the NLS equation}
For the nonlinear Schrödinger equation, the general functions in \eqref{eq:quasigeneral} become
\(\mathcal{L}q = -\ii q_{xx}\)
and
\(\mathcal{N}(q,q_x) = -2\ii |q|^2q.\)
Therefore, \(M = -\ii D_J^2\) and the system of ordinary differential equations in \eqref{eq:genode} becomes
\[\vec{a}'(t)=\ee^{-\ii tD_J^2}F(\ee^{\ii tD_J^2}\vec{a}(t)),\]
for \(\vec{a}(t)=\ee^{-\ii tD_J^2}\vec{c}(t)\) where
\[F(\vec{c})=2\ii\mathcal{F}_J(\mathcal{F}_J^{-1}(\vec{c})\cdot|\mathcal{F}_J^{-1}(\vec{c})|^2)\]
and, again, $\cdot$ denotes the entrywise multiplication of two vectors.

\subsection{Artificially-damped solutions of the NLS equation}

We compute solutions to \eqref{eq:nls} at \(t = 150\) that are accurate (errors on the order of \(10^{-8}\)) on the interval \([-99.85,100.05]\), using a time step of 0.01 for all computations. We effectively solve the modified equation,
\begin{equation}
    \ii q_t+q_{xx}+2|q|^2q=-k_2(1-\gamma(x))q,
    \label{eq:modnls}
\end{equation}
by multiplying the solution values by an even \(\gamma(x)\), defined as
\begin{equation}
    \gamma(x) = 1 - (\sigma(x) + \sigma(-x)),
    \label{eq:newgamma}
\end{equation}
where \(\sigma(x)\) is the same as in \eqref{eq:sigma}.
\begin{figure}
    \centering
    \includegraphics[scale = 0.4]{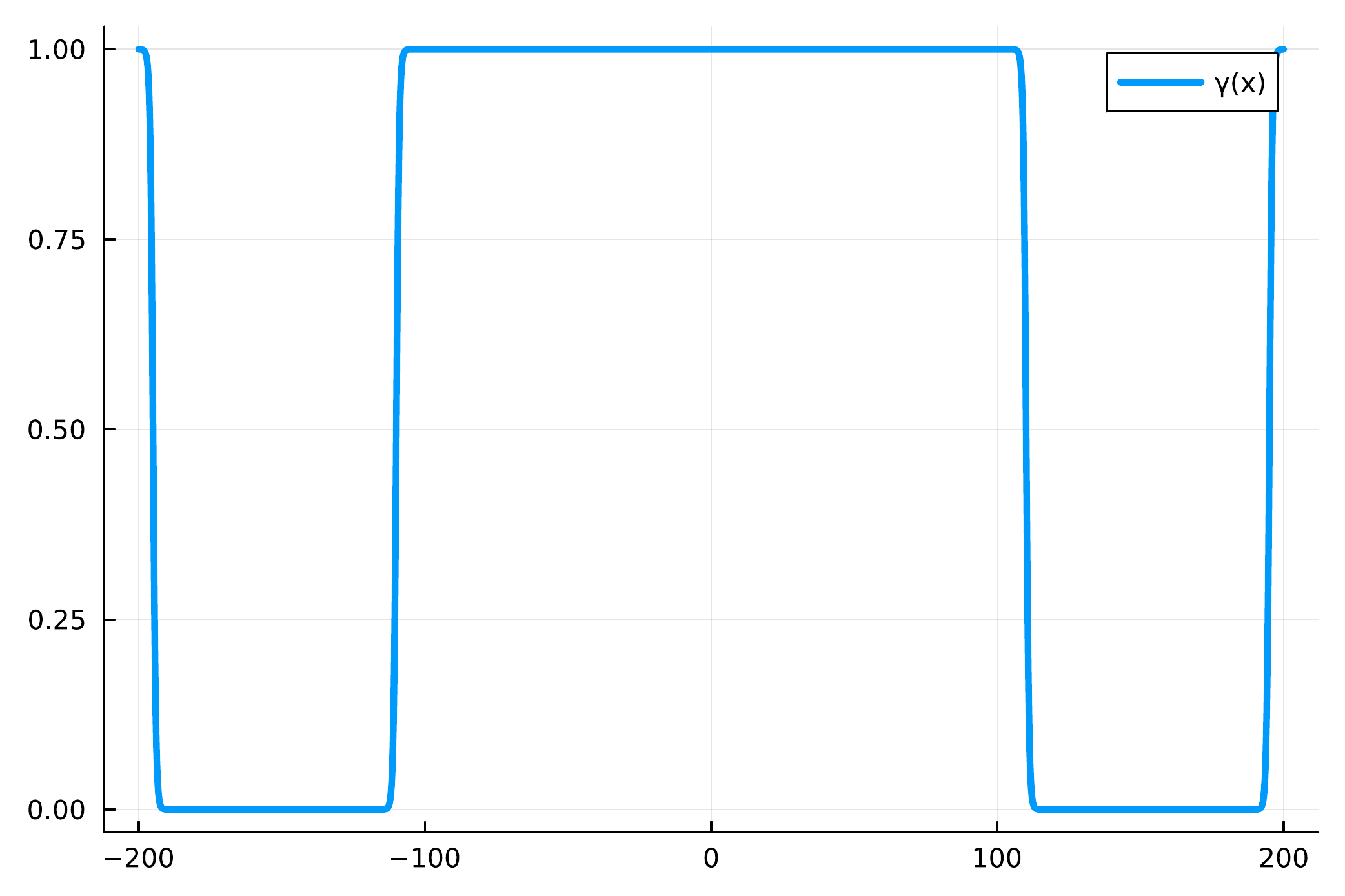}
    \caption{Even damping factor, \(\gamma(x)\).}
    \label{fig:fig8}
\end{figure}
We chose an even \(\gamma(x)\) function, shown in \eqref{fig:fig8}, due to the more symmetric nature of the real part of the solution. As in the previous section, we perform this multiplication every 1000 time steps. We do not include the \(k_1(\sigma(x)q_{x})_x\) term in the modified equation because of the increase in runtime that using Strang-splitting  to solve the heat equation causes. For the KdV equation, we found that using this form of damping significantly improved the accuracy of solutions on smaller intervals, but for the NLS equation, we were able to achieve the desired accuracy without this term.

\begin{figure}
    \centering
    \includegraphics[width = 0.8\linewidth]{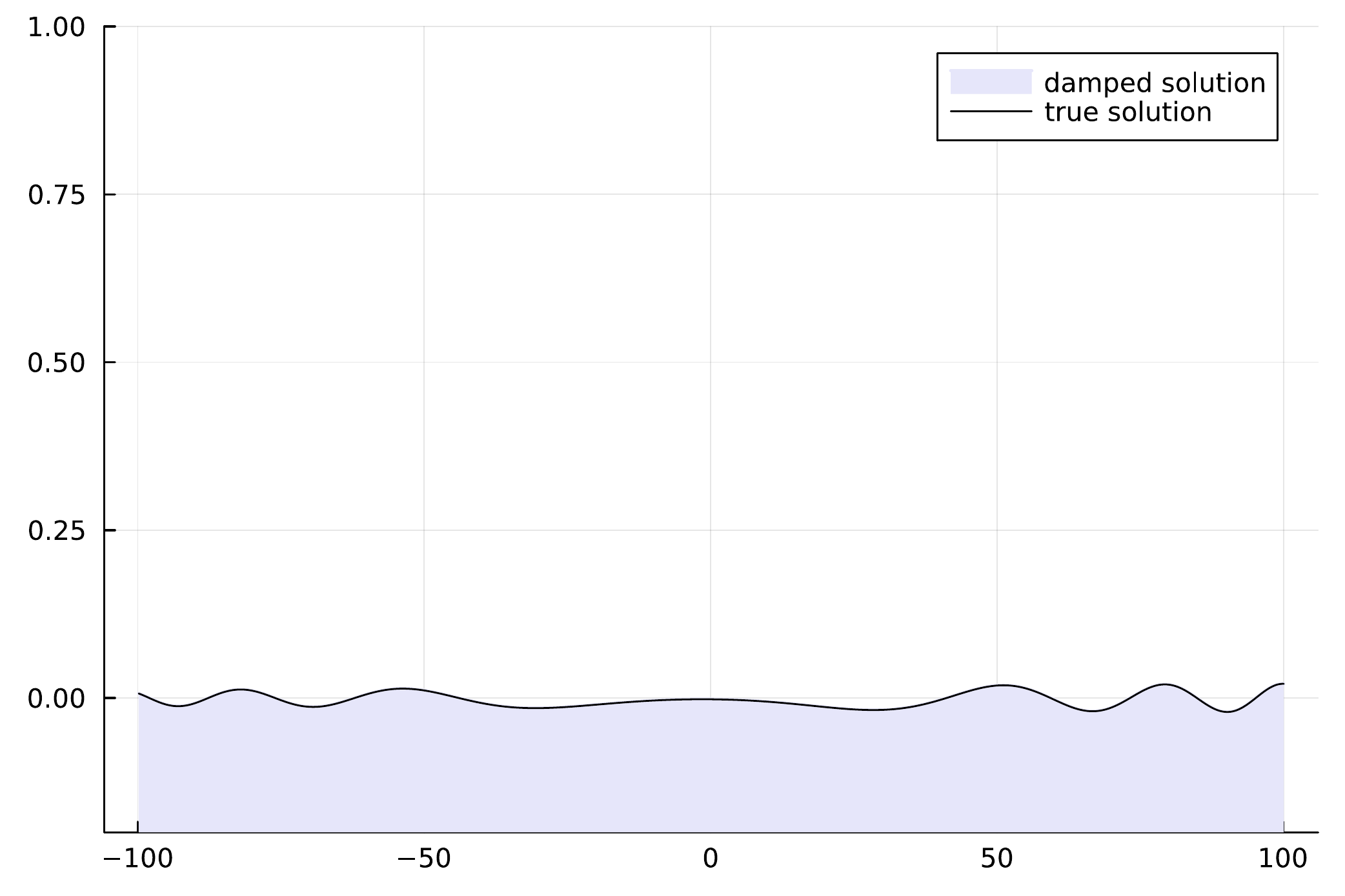}
    \caption{Damped real part of the solution to the NLS equation computed on \([-1200,1200]\).}
    \label{fig:fig9}
\end{figure}

The real part of the solution shown in Figure~\ref{fig:fig9} is computed on \([-1200, 1200]\) with \(2^{13}\) grid points. This takes 50 seconds and has a maximum error\footnote{The true solution in this section comes from \cite{trogdonschrodinger}.} of \(10^{-8}\). Without damping, an interval of \([-2500,2500]\) and \(2^{14}\) grid points is required to achieve the same order of accuracy. This solution takes 91 seconds to compute, increasing the runtime by \(82\%\) in comparison to the damped solution. Table \ref{tab:nlstable} contains further results comparing the errors and runtimes of damped and undamped solutions computed on varying intervals. Again, we see that the error saturates on the order of \(10^{-8}\) due to our use of fourth-order Runge-Kutta and a time step of 0.01.

\begin{table}[h!]
      \centering
      \begin{tabular}{l r r c}
        \toprule
        \textbf{Parameters} & \textbf{Undamped Error} & \textbf{Damped Error}\\
        \midrule
        \(L = 100,\;m = 2^9\) & \(0.6\) (5 sec) & \(0.02\) (5 sec) \\
        \(L = 200,\;m = 2^{10}\) &  0.5 (8 sec) &  \(0.005\) (8 sec)\\
        \(L = 600,\; m = 2^{12}\) &  0.002 (26 sec) & \(3 \cdot 10^{-5}\) (27 sec)\\
        \(L = 1200,\; m = 2^{13}\) &  \(1\cdot 10^{-5}\) (46 sec) & \(1 \cdot 10^{-8}\) (50 sec)\\
        \(L = 2500,\; m = 2^{14}\) &  \(1\cdot 10^{-8}\) (91 sec) &  \(5 \cdot 10^{-8}\) (95 sec)\\
        \bottomrule
      \end{tabular}
      \caption{Solutions evaluated on an interval \([-L,L]\) with m grid points.}
      \label{tab:nlstable}
    \end{table}

Now having shown that our damping technique can improve the runtime of computing accurate solutions to the Korteweg-de Vries equation and the nonlinear Schrödinger equation, we turn our attention to a slightly more complicated problem.

\section{The Riemann problem for the KdV equation}
\label{sec:shock}

We consider the Riemann problem for the KdV equation \eqref{eq:kdv}, 
\begin{equation}
    q_t+qq_x+\epsilon^2q_{xxx} = 0,
    \label{eq:cauchy}
\end{equation}
with \(\epsilon=10^{-1.5}\) and the initial condition, \(q(x,0)=\frac{1}{1+\ee^{10x}}\). A similar problem is discussed in \cite{cauchyproblem}. Since the Fourier method requires periodicity and \(q(x,0)\) is not periodic, as described in \cite{Sprenger2017} and \cite{hoefershearer}, we use the fact that the derivative of \(q(x,0)\) is nearly periodic and obtain the Fourier coefficients \(\vec{c}(t)\) of \(u(x,t)=q_x(x,t)\), which satisfies the following nonlinear PDE,
\begin{equation}
    u_t+(\partial_x^{-1}u)u_x+u^2+\epsilon^2u_{xxx}=0.
    \label{eq:periodic}
\end{equation}

\subsection{The Fourier method for the Riemann problem}
Applying the Fourier method in Figure~\ref{sec:fourier} to the problem, we have that the general functions in \eqref{eq:quasigeneral} are
\(\mathcal{L}(q) = \epsilon^2q_{xxx}\)
and
\(\mathcal{N}(q,q_x) = qq_x.\)
Therefore, \(M = \epsilon^2D_J^3\) and the system of ordinary differential equations in \eqref{eq:genode} becomes
\[\vec{a}'(t)=\ee^{\epsilon^2D_J^3t}F(\ee^{-\epsilon^2D_J^3t}\vec{a}(t))\]
with \(\vec{a}(t)=\ee^{\epsilon^2D_J^3t}\vec{c}(t)\) and 
\[F(\vec{c})=-\mathcal{F}_J\left(\mathcal{F}_J^{-1}(D_J\cdot \vec{c})\cdot H(\vec{c}) + \mathcal{F}_J^{-1}(\vec{c})^2\right),\]
where \(H(\vec{c})\) computes the antiderivative of \(u(\vec{x},t)\).

\subsection{Defining \(H(\vec c)\)}
We begin by expressing \(q(x,t)\) as the integral of the approximate Fourier series of \(u(x,t)\).
\[q(x,t)\approx\int_{-L}^x\sum_{n=-N}^Nc_n(t)\ee^{\frac{\ii n\pi}{L}s}\dd s+C_-,\]
We normalize the antiderivative of \(u(x,t)\) at \(-L\) by defining the integration constant as \(C_-=q(-L,0)\). Pulling out the \(n=0\) term from the sum and integrating it separately yields
\begin{align*}
q(x,t)&\approx C_-+c_0(t)(x+L)+\sum_{n=-N,n\neq0}^N\left(\frac{L}{\ii n\pi}\right)c_n(t)\ee^{\frac{\ii n\pi}{L}s}\biggr|_{-L}^x\\
&=c_0(t)(x+L)+C_--\sum_{n=-N,n\neq0}^N\left(\frac{L}{\ii n\pi}\right)c_n(t)\ee^{-\ii n\pi}+\sum_{n=-N,n\neq0}^N\left(\frac{L}{\ii n\pi}\right)c_n(t)\ee^{\frac{\ii n\pi}{L}x}\\
&=c_0(t)(x+L)+\sum_{n=-N}^N S_n(t)\ee^{\frac{\ii n\pi}{L}x}
\end{align*}
where \(S_n(t)\) is defined as
\[S_n(t) =  \begin{cases} 
      \frac{L}{\ii n\pi}c_n(t) & n\neq0, \\
       C_--\sum_{n=-N,n\neq0}^N\left(\frac{L}{\ii n\pi}\right)c_n(t)\ee^{-\ii n\pi} & n=0.
   \end{cases}
\] 
We can then define the function \(H(\vec{c})\) through Algorithm~\ref{alg:integrate}, which takes in a vector of the Fourier coefficients of \(u(x,t)\) and computes a vector of values of \(q(x,t)\) along the interval \([-L,L]\) with \(m\) evenly spaced grid points. In this algorithm we use the notation that for a diagonal matrix $D$, $D^\dagger$ is defined by $(D^\dagger)_{ij} = D_{ij}^{-1}$ if $i =j$, $D_{ii} \neq 0$ and $(D^\dagger)_{ij} = 0$ otherwise.
\begin{algorithm}
\caption{Integrate}
\label{alg:integrate}
\begin{algorithmic}[1]
    \Function{H}{$\vec{c}$}
        \State{$c_0 =$ Fourier coefficient corresponding to $j = 0$}
        \State{$Q := D_J^{\dagger}$}
        \State{$S=Q\vec{c}$}
        \State{$\vec{q}=\frac{c_0}{m}(x+L)+\mathcal{F}_J^{-1}(S)$}
        \State{$\vec{q}=\vec{q}-\vec{q}_{-L}+C_-$}
        \State{\textbf{return }$\vec{q}$}
    \EndFunction
\end{algorithmic}
\end{algorithm}

\subsection{Solutions to the Riemann problem}

\begin{figure}[htp]
    \centering
    \includegraphics[width=0.9\linewidth]{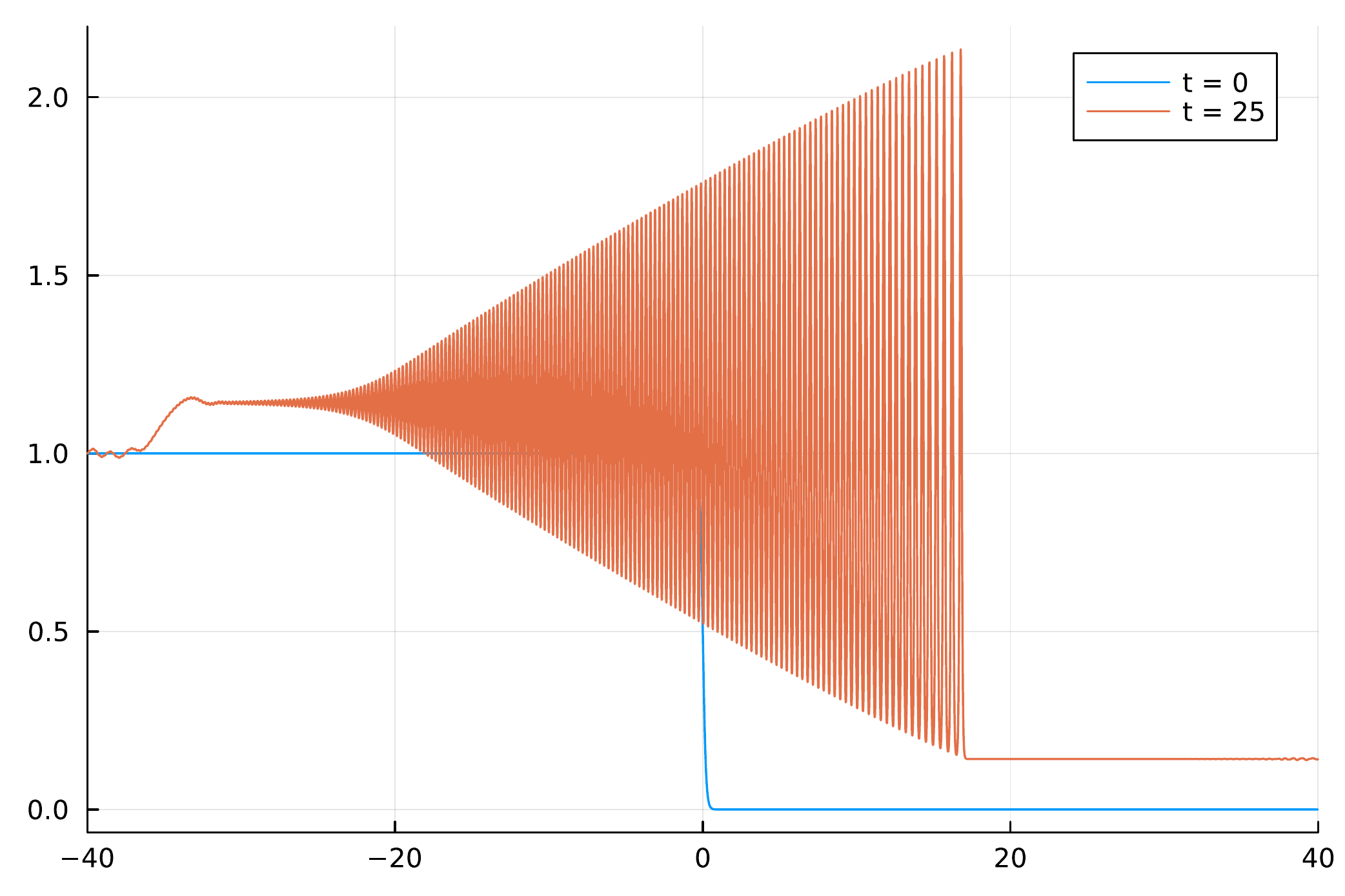}
    \caption{Undamped solution of the Riemann problem problem computed on \([-40,40]\).}
    \label{fig:shock1}
\end{figure}

\begin{figure}[htp]
    \centering
    \includegraphics[width=0.9\linewidth]{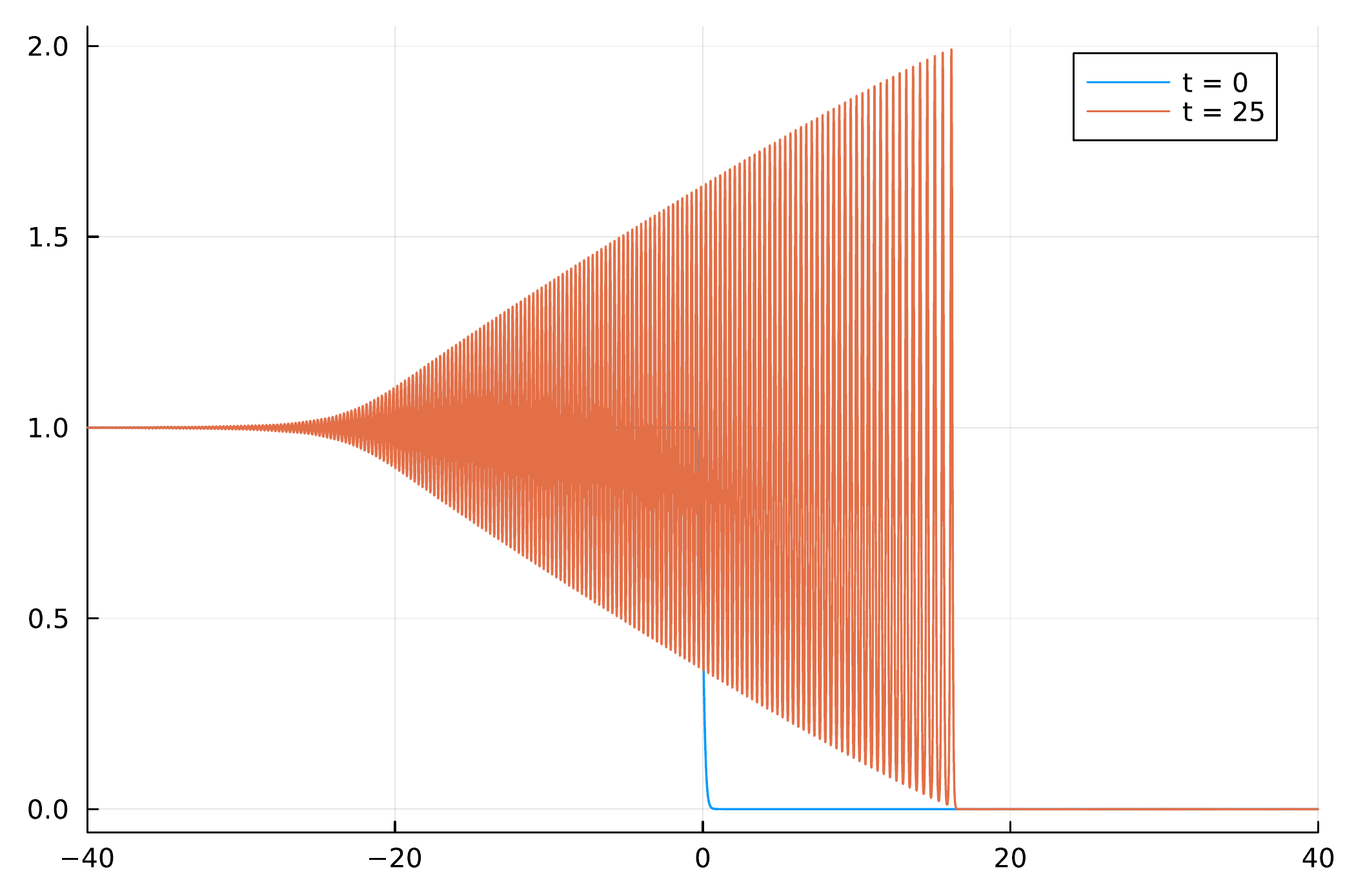}
    \caption{Undamped solution of the Riemann problem computed on \([-60,60]\).}
    \label{fig:fig6}
\end{figure}

\begin{figure}[htp]
    \centering
    \includegraphics[width=0.9\linewidth]{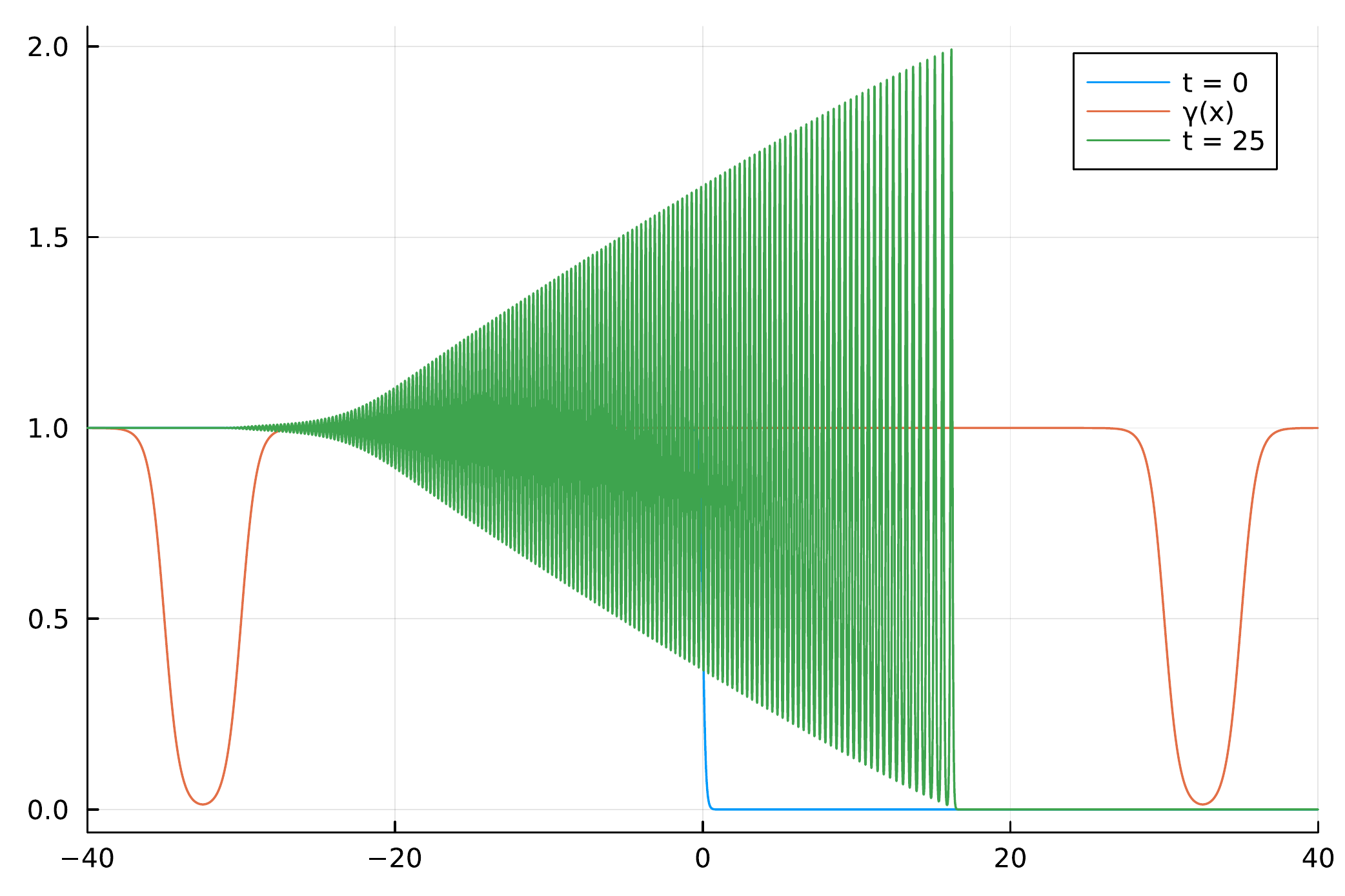}
    \caption{Damped solution of the Riemann problem computed on \([-40,40]\).}
    \label{fig:fig7}
\end{figure}

Unlike in previous problems, we do not have a ground-truth solution to compare our computations to because the numerical inverse scattering transform for this problem currently unable to evaluate solutions \cite{Bilman2020} at a sufficiently large time $t$. But we do know what solutions should look like based on \cite{cauchyproblem}, for example. We begin by computing a solution\footnote{All solutions in this section are computed with a time step of 0.01.} to \eqref{eq:cauchy} at \(t = 25\) on the interval \([-40,40]\) with \(2^{12}\) grid points. A higher density of grid points than in previous sections is required due to the solution forming a shock wave. Without any damping, this computation takes 68 seconds and produces the plot in Figure~\ref{fig:shock1}. 
Wrap-around oscillations have caused the entire shock wave to shift upwards unnaturally. We can improve the solution quality by increasing the computational domain to \([-60,60]\), shown in Figure~\ref{fig:fig6}, or introducing artificial damping, shown in Figure~\ref{fig:fig7}. Following Section~\ref{sec:damping1}, we damp the shock wave at the edges of the domain by effectively solving
\begin{equation*}
    u_t+(\partial_x^{-1}u)u_x+u^2+\epsilon^2u_{xxx}=-k_2(1-\gamma(x))u
\end{equation*}
and defining \(\gamma(x)\) as in \eqref{eq:newgamma}. Again, we choose to set \(k_1 = 0\) and do not use Strang-splitting to solve the heat equation because we found it possible to get a qualitatively accurate solution at a faster time without it. Both methods solve the problem, but increasing the computational domain nearly doubles the runtime to 120 seconds, while using damping maintains a runtime of 68 seconds.


\begin{figure}[htp]
    \centering
    \includegraphics[width=0.9\linewidth]{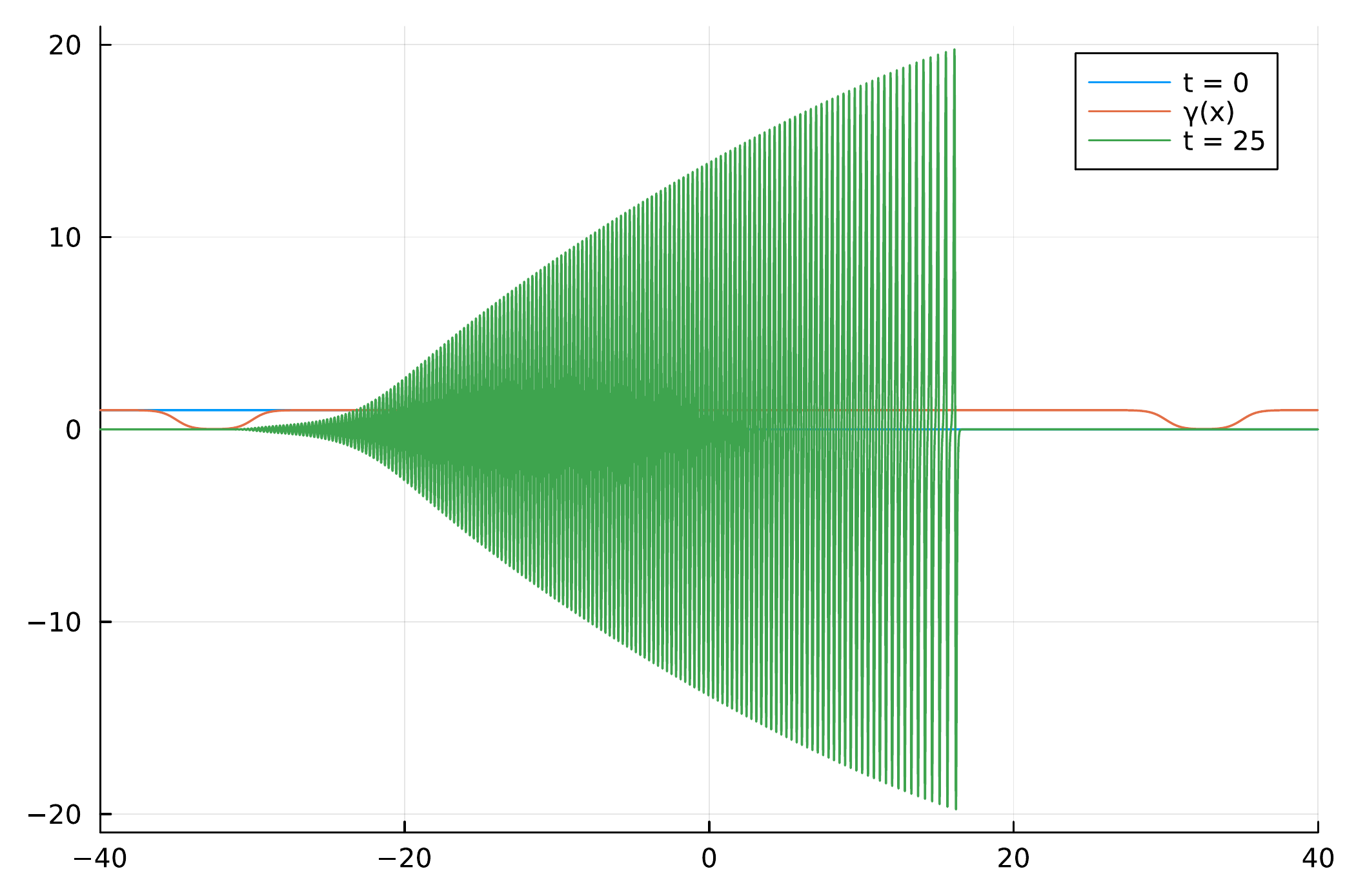}
    \caption{Derivative of the damped solution of the Riemann problem computed on \([-40,40]\).}
    \label{fig:shockderiv}
\end{figure}

Our technique relies on integrating the derivative of the solution at each time step. In Figure~\ref{fig:shockderiv}, we plot the derivative of the solution at \(t = 25\). The derivative looks like the solution scaled up by 10 - a result of the chosen initial condition, \(q(x,0)=\frac{1}{1+\ee^{10x}}\). We also observe that despite the high frequency and amplitude of the derivative, our integrating function, \(H(\vec c)\), produces a qualitatively accurate shock wave solution.

\section{The Kawahara equation}
\label{sec:kawa}
Finally, we use our damping methods to solve the Kawahara equation \cite{ogkawahara},
\begin{equation}
    q_t+qq_x+q_{xxx}+q_{xxxxx}=0,
    \label{eq:kawa}
\end{equation}
with the initial condition \(q(x,0)=\frac{1}{1+\ee^{10x}}-1\). Again, we have an aperiodic initial condition so we rewrite \eqref{eq:kawa} to solve for \(u(x,t) = q_x(x,t)\),
\begin{equation}
    u_t+(\partial_x^{-1}u)u_x+u^2+u_{xxx}+u_{xxxxx}=0.
    \label{eq:derivkawa}
\end{equation}

\subsection{The Fourier method for the Kawahara equation}
In the notation of Section~\ref{sec:fourier}, we define \(\mathcal{L}q = q_{xxx}+q_{xxxxx}\) and \(\mathcal{N}(q,q_x) = qq_x\) such that the Kawahara equation is in the form of \eqref{eq:quasigeneral}. From this, we have \(M = D_J^3+D_J^5\) and the system of ordinary differential equations in \eqref{eq:genode} becomes an approximate finite-dimensional ODE system,
\[\vec{a}'(t)=\ee^{(D_J^3+D_J^5)t}F(\ee^{-(D_J^3+D_J^5)t}\vec{a}(t)),\]
where \(\vec{a}(t) = \ee^{(D_J^3+D_J^5)t}\vec{c}(t)\) and 
\[F(\vec{c})=-\mathcal{F}_J\left(\mathcal{F}_J^{-1}(D_J\cdot \vec{c})\cdot H(\vec{c}) + \mathcal{F}_J^{-1}(\vec{c})^2\right).\]
\(H(\vec{c})\) is still defined through Algorithm~\ref{alg:integrate}.

\subsection{Solutions to the Kawahara equation}

Sprenger and Hoefer solve the Kawahara equation \eqref{eq:kawa} in \cite{Sprenger2017} with the initial condition,
\[q(x,0)= \begin{cases} 
    0 & x<0, \\
    -1 & x\geq 0,
    \end{cases}
\]
and the same method described in Section~\ref{sec:shock}. Motivated by their work, we compute a reference solution at \(t = 24\) on \([-10^5,10^5]\). We use a time step of 0.0005 and \(m=2^{22}\) grid points; we approximate their discontinuous initial condition with \(q(x,0)=\frac{1}{1+\ee^{10x}}-1\). This computation took over 30 hours\footnote{This computation was done on an Intel® Core™ i7-6700 Processor with 32 GB RAM, a 3.40GHz CPU, and the CentOS operating system.}. We now compute solutions with and without damping on \([-1000,1000]\) with \(2^{15}\) grid points. We introduce damping by solving a modified version of \eqref{eq:derivkawa},
\begin{equation}
    u_t+(\partial_x^{-1}u)u_x+u^2+u_{xxx}+u_{xxxxx}=-k_2(1-\gamma(x))u,
    \label{eq:derivkawadamped}
\end{equation}
where \(\gamma(x)\) is still defined to be an even function as in \eqref{eq:newgamma}. We multiply the solution values by \(\gamma(x)\) every 100 time steps (\(f_2 = 100\)). 
\begin{figure}
    \centering
    \includegraphics[scale = 0.55]{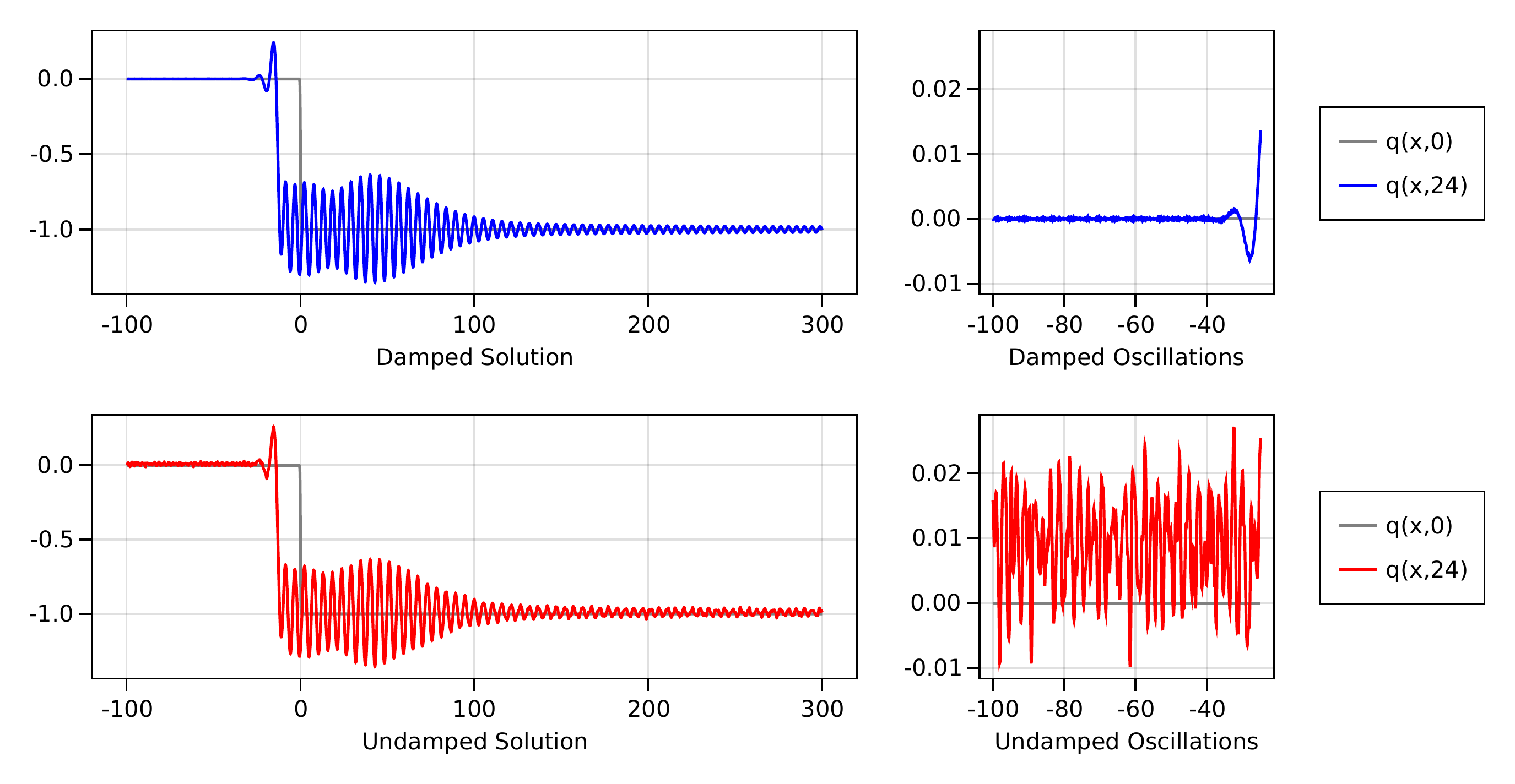}
    \caption{Solutions to the Kawahara equation at \(t = 24\) computed on \([-1000,1000]\).}
    \label{fig:kawacomparison}
\end{figure}
More aggressive damping is required for this problem due to the fifth-order term that causes higher velocity dispersion. We compare the damped and undamped solutions in Figure~\ref{fig:kawacomparison}. The plots on the right show a smaller subset of the interval to highlight the effect of damping on unwanted oscillations. The damped solution took 444 seconds to compute and has a maximum error of 0.001 from the reference solution. The undamped solution had a shorter runtime of 492 seconds, but produced a maximum error of 0.03.

\section{Summary of results}

We have developed a general method for modeling nonlinear dispersive PDEs on the line with a finite computational domain. Our technique utilizes the Fourier method which is a very efficient numerical tool; however, it forces periodicity onto the solution. In order to accurately model the dispersion of solutions, we modify the PDE to include damping terms, \(k_1(\sigma(x)q_x)_x\) and \(k_2(1-\gamma(x))q)\). The former solves the heat equation in the desired damping region, specified by the choice of \(\sigma(x)\), to diffuse the solution near the edges of the interval. The latter causes rapid exponential decay in the region determined by \(\gamma(x)\). We found that \(\sigma(x)\) damping works particularly well for dispersive tails, while \(\gamma(x)\) damping is better suited to damping solitons. Different choices of damping parameters will vary the accuracy and computation of approximate solutions. 

% \section{Conclusions}
% We demonstrated the effectiveness of our method on four dispersive PDEs. First, we solved the KdV equation, using heat equation damping for the dispersive tail and exponential decay damping for the soliton. Next, we solved the NLS equation, using only \(\gamma(x)\) damping on both sides of the computational domain. We then solved a Riemann problem for the KdV equation with a step-like initial condition that produced a shock wave solution. Again, we found that \(\gamma(x)\) damping alone was enough to get a qualitatively accurate solution at an improved runtime. Finally, we applied our method to a Riemann problem for the Kawahara equation. By applying our damping technique to four different PDEs, we have highlighted its generality and demonstrated how it can be modified to suit a specific problem.

\appendix
\section{Further applications} 

\begin{figure}[H]
\begin{subfigure}{.5\textwidth}
  \centering
  \includegraphics[width=.7\linewidth]{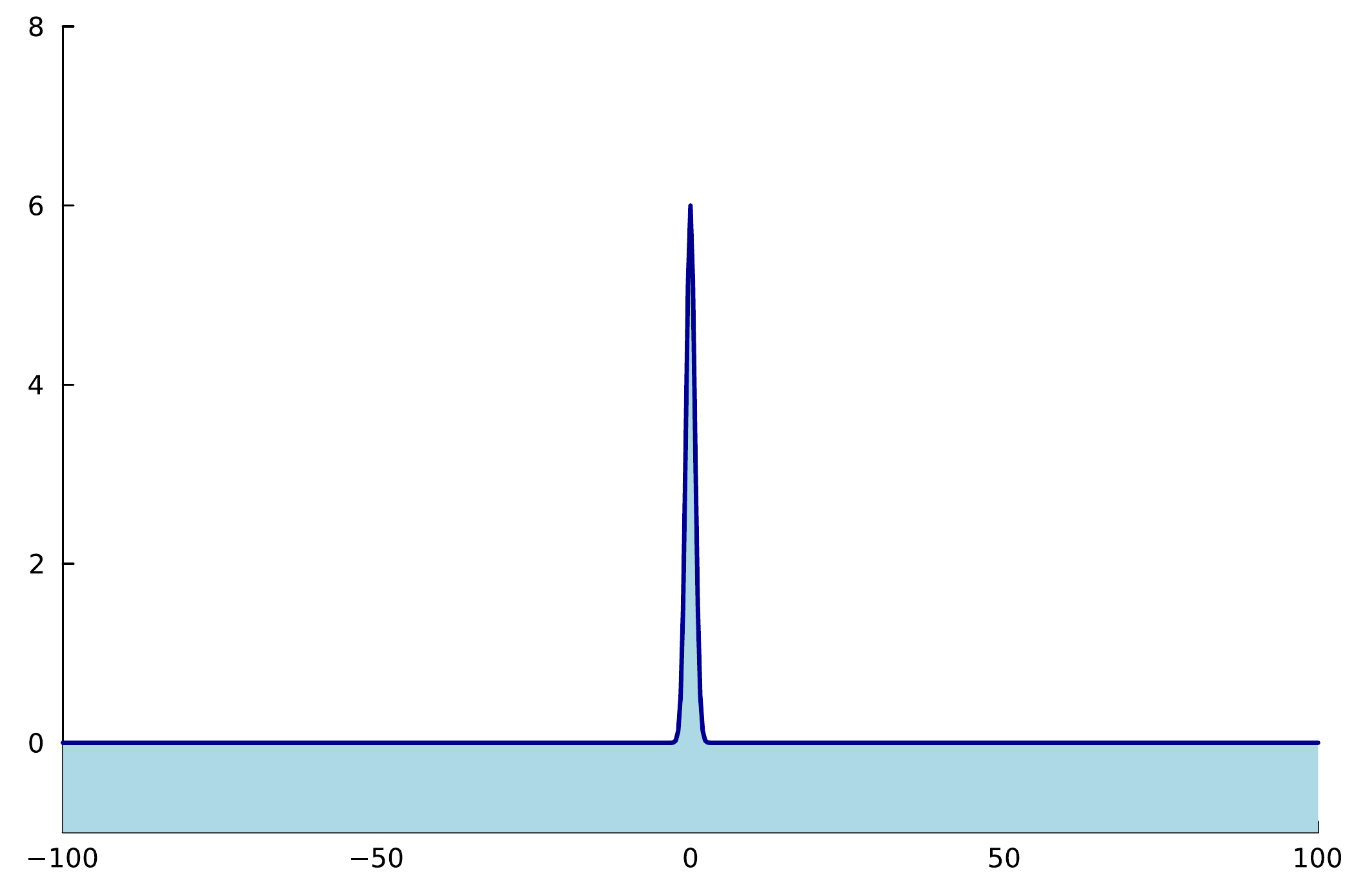}
  \caption{Initial condition, \(q(x,0) = 6\ee^{-x^2}\).}
  \label{fig:sfig1}
\end{subfigure}%
\begin{subfigure}{.5\textwidth}
  \centering
  \includegraphics[width=.7\linewidth]{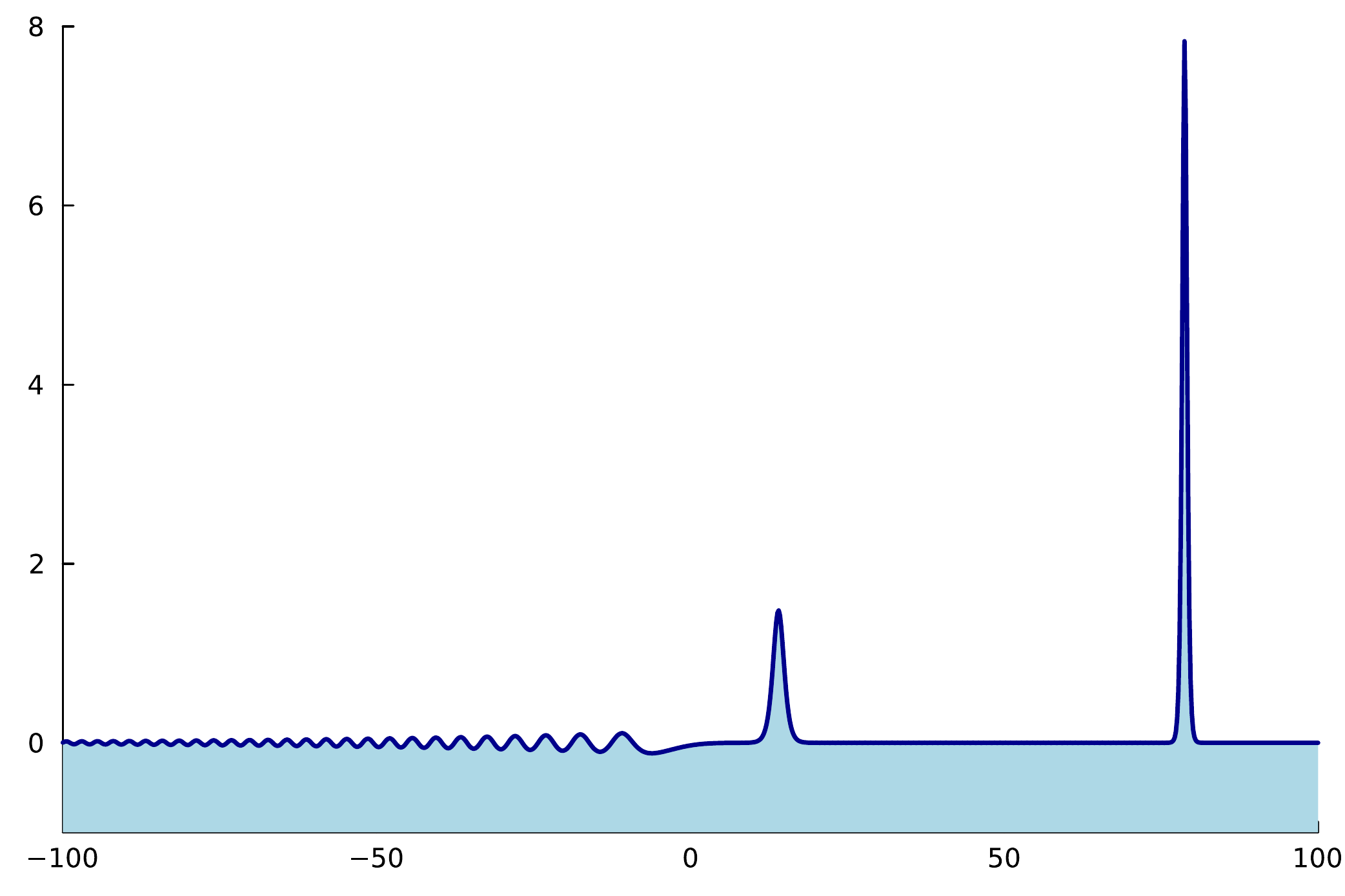}
  \caption{Solution at \(t = 5\).}
  \label{fig:sfig2}
\end{subfigure}
\caption{Solution to the KdV equation \eqref{eq:kdv}. \textbf{Parameters:} \(L = 200\), \(m = 2^{11}\), \(\triangle t = 0.001\), \(\sigma(x) = \) \eqref{eq:sigma}, \(\gamma(x) = \) \eqref{eq:gamma}, \(k_1 = 1\), \(f_1 = 1\), \(f_2 = 1000\). \textbf{Runtime:} 28 seconds.}
\label{fig:twosoliton}
\end{figure}

\begin{figure}[H]
\begin{subfigure}{.5\textwidth}
  \centering
  \includegraphics[width=.7\linewidth]{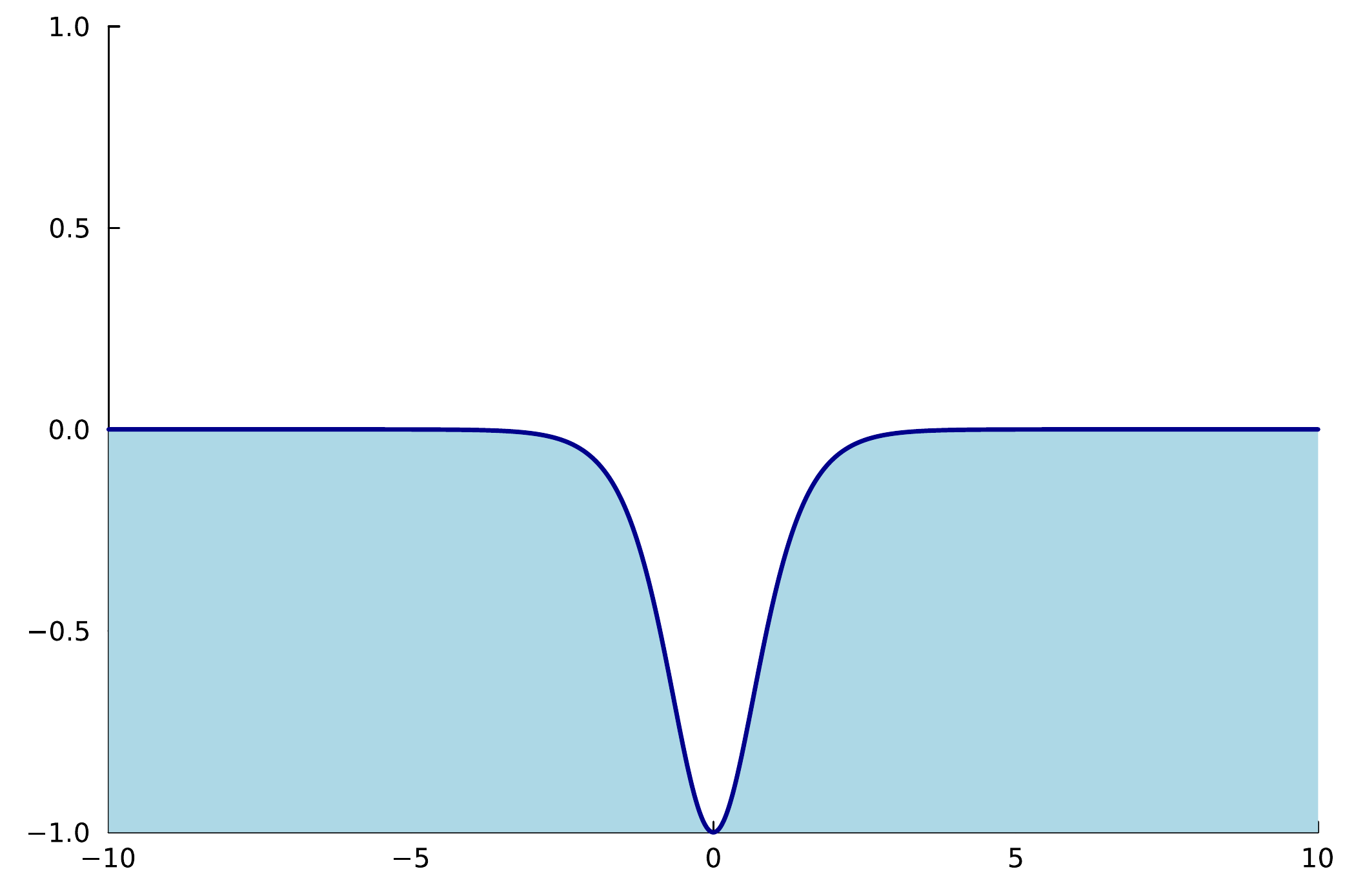}
  \caption{Initial condition, \(q(x,0) = -\frac{1}{\cosh^2(x)}\) (as in \cite{cauchyproblem}).}
  \label{fig:sfig14}
\end{subfigure}%
\begin{subfigure}{.5\textwidth}
  \centering
  \includegraphics[width=.7\linewidth]{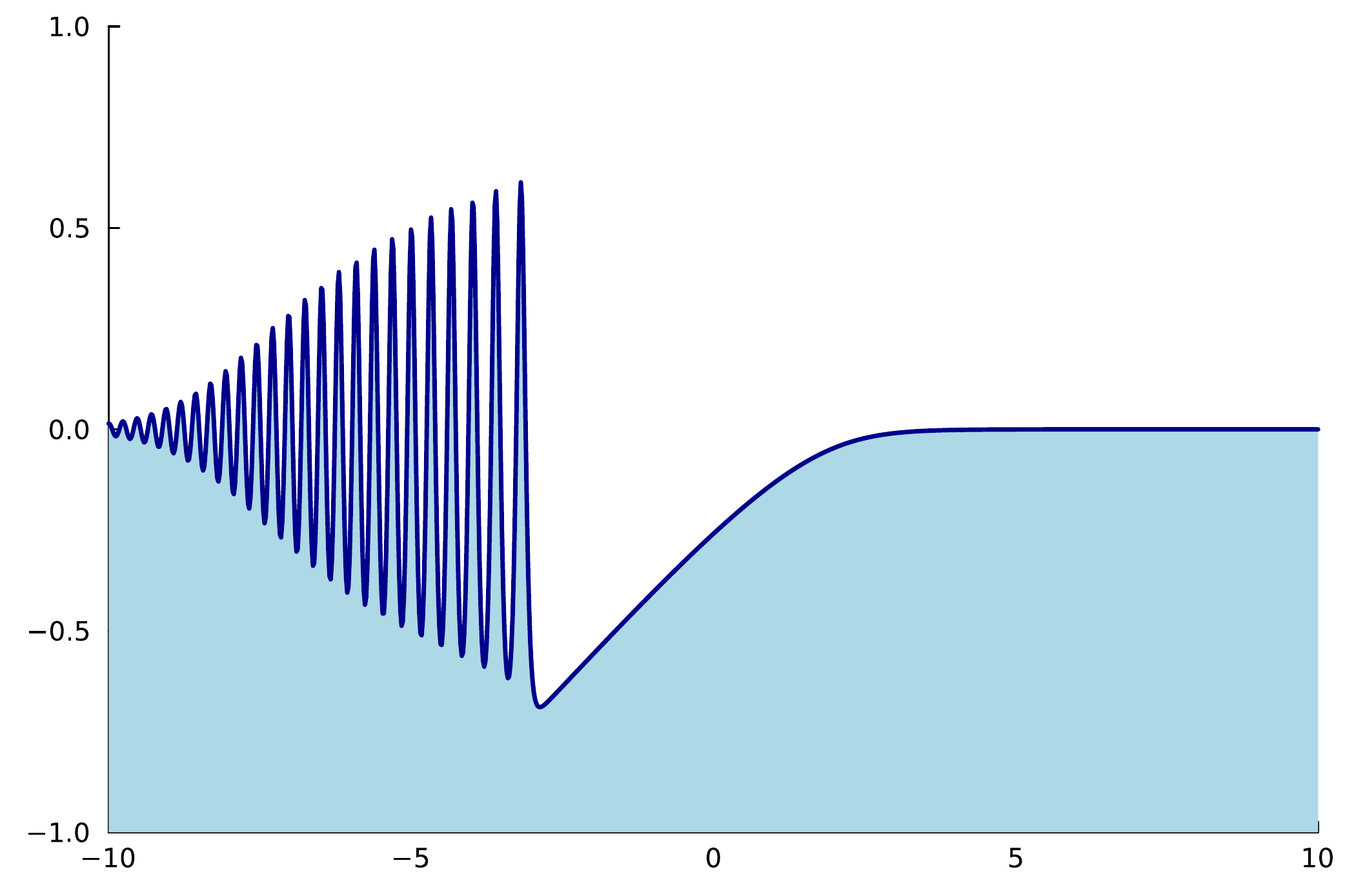}
  \caption{Solution at \(t = 5\).}
  \label{fig:sfig23}
\end{subfigure}
\caption{Solution to the Riemann problem for the KdV equation \eqref{eq:cauchy}. \textbf{Parameters:} \(L = 40\), \(m = 2^{12}\), \(\triangle t = 0.001\), \(\sigma(x) = \) \eqref{eq:sigma}, \(\gamma(x) = \) \eqref{eq:newgamma}, \(k_1 = 0\), \(f_2 = 1000\). \textbf{Runtime:} 15 seconds.}
\label{fig:riemannappend}
\end{figure}

\begin{figure}[H]
\begin{subfigure}{.5\textwidth}
  \centering
  \includegraphics[width=.7\linewidth]{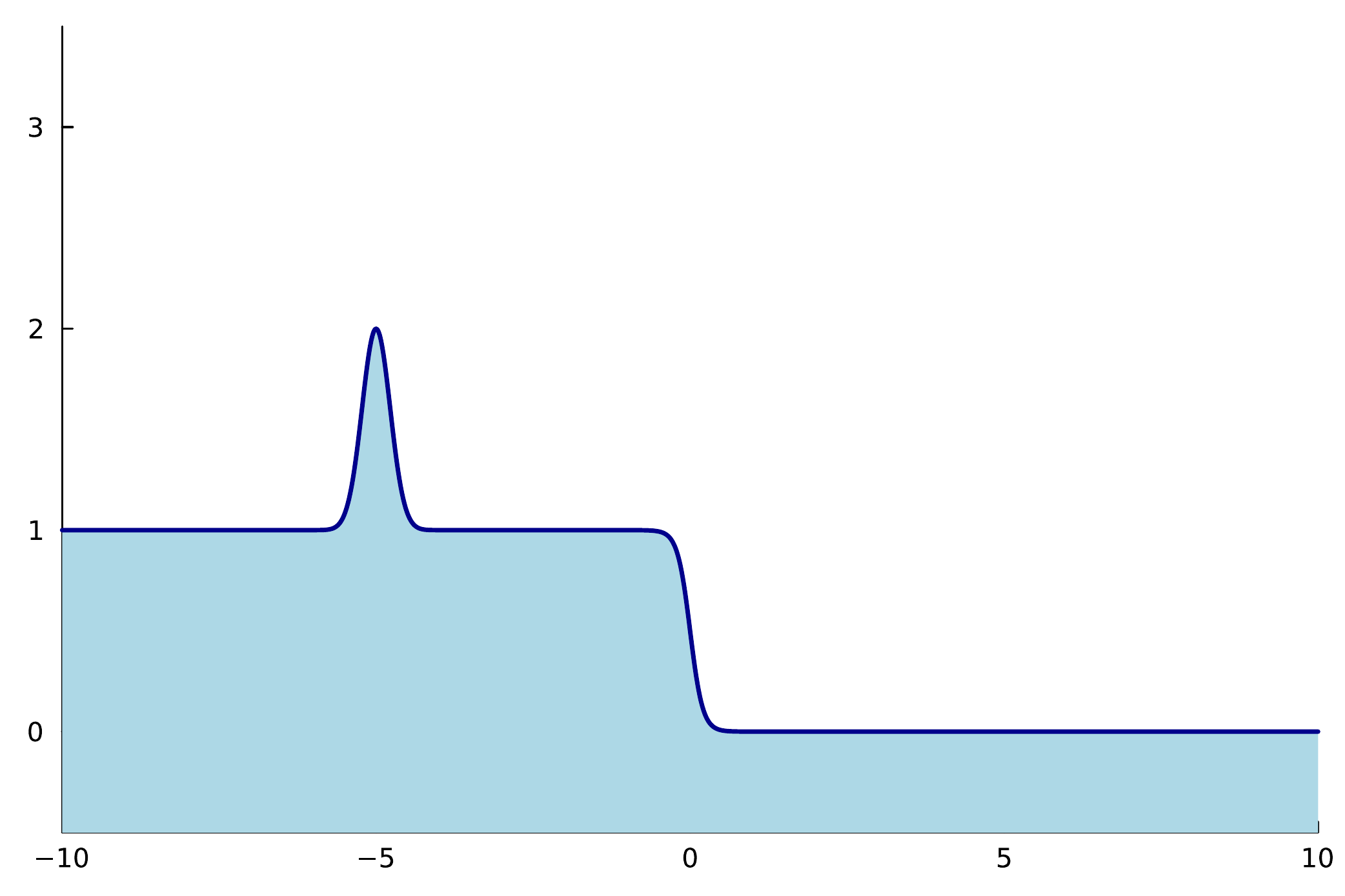}
  \caption{Initial condition, \(q(x,0) = -\frac{1}{1+\ee^{-10x}}+1+\ee^{-10(x+5)^2}\).}
  \label{fig:sfig18}
\end{subfigure}%
\begin{subfigure}{.5\textwidth}
  \centering
  \includegraphics[width=.7\linewidth]{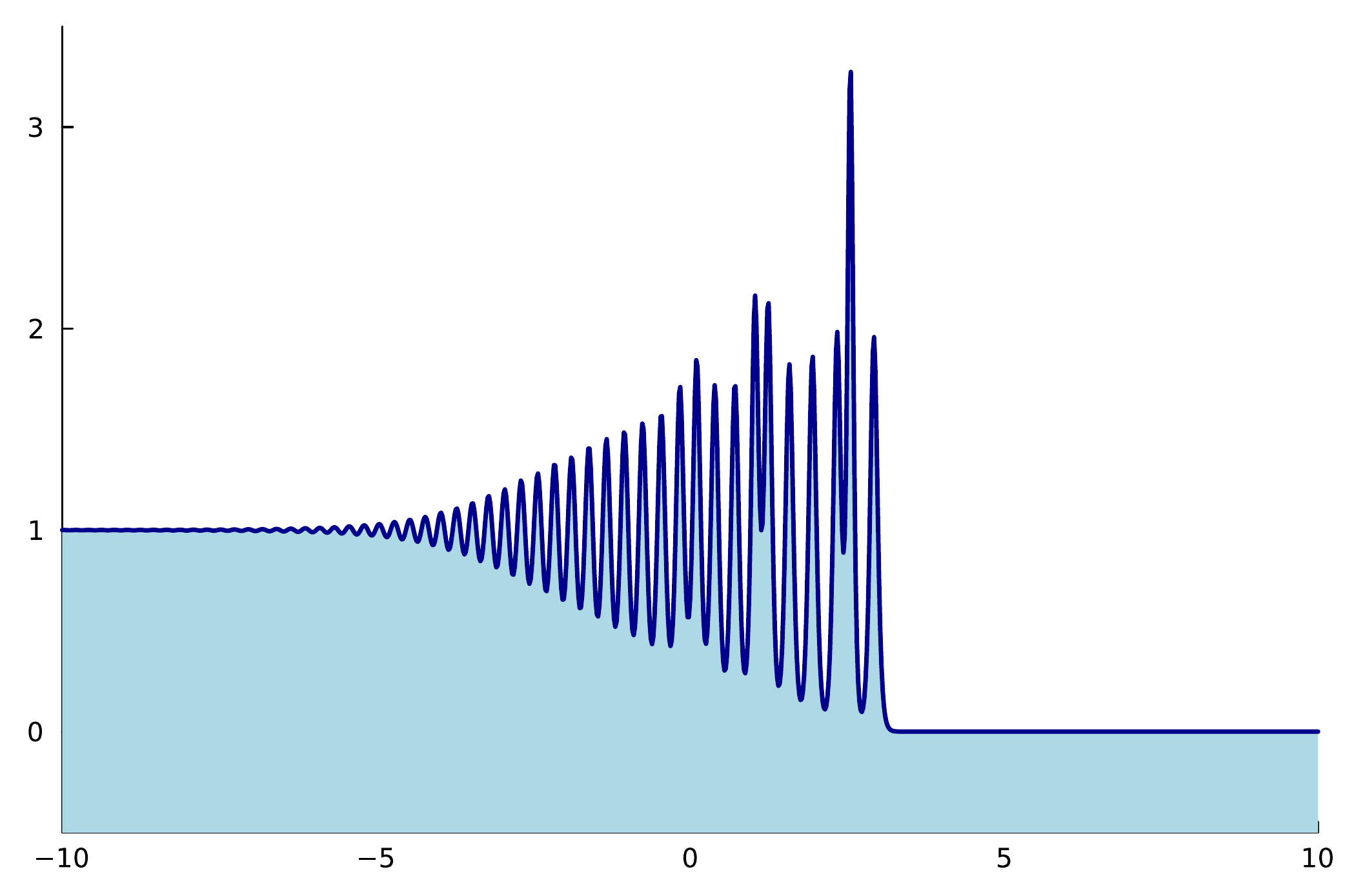}
  \caption{Solution at \(t = 5\).}
  \label{fig:sfig27}
\end{subfigure}
\caption{Solution to the Riemann problem for the KdV equation \eqref{eq:cauchy}. \textbf{Parameters:} \(L = 40\), \(m = 2^{12}\), \(\triangle t = 0.001\), \(\sigma(x) = \) \eqref{eq:sigma}, \(\gamma(x) = \) \eqref{eq:newgamma}, \(k_1 = 0\), \(f_2 = 1000\). \textbf{Runtime:} 18 seconds.}
\label{fig:solitonshelf}
\end{figure}

\begin{figure}[H]
\begin{subfigure}{.5\textwidth}
  \centering
  \includegraphics[width=.7\linewidth]{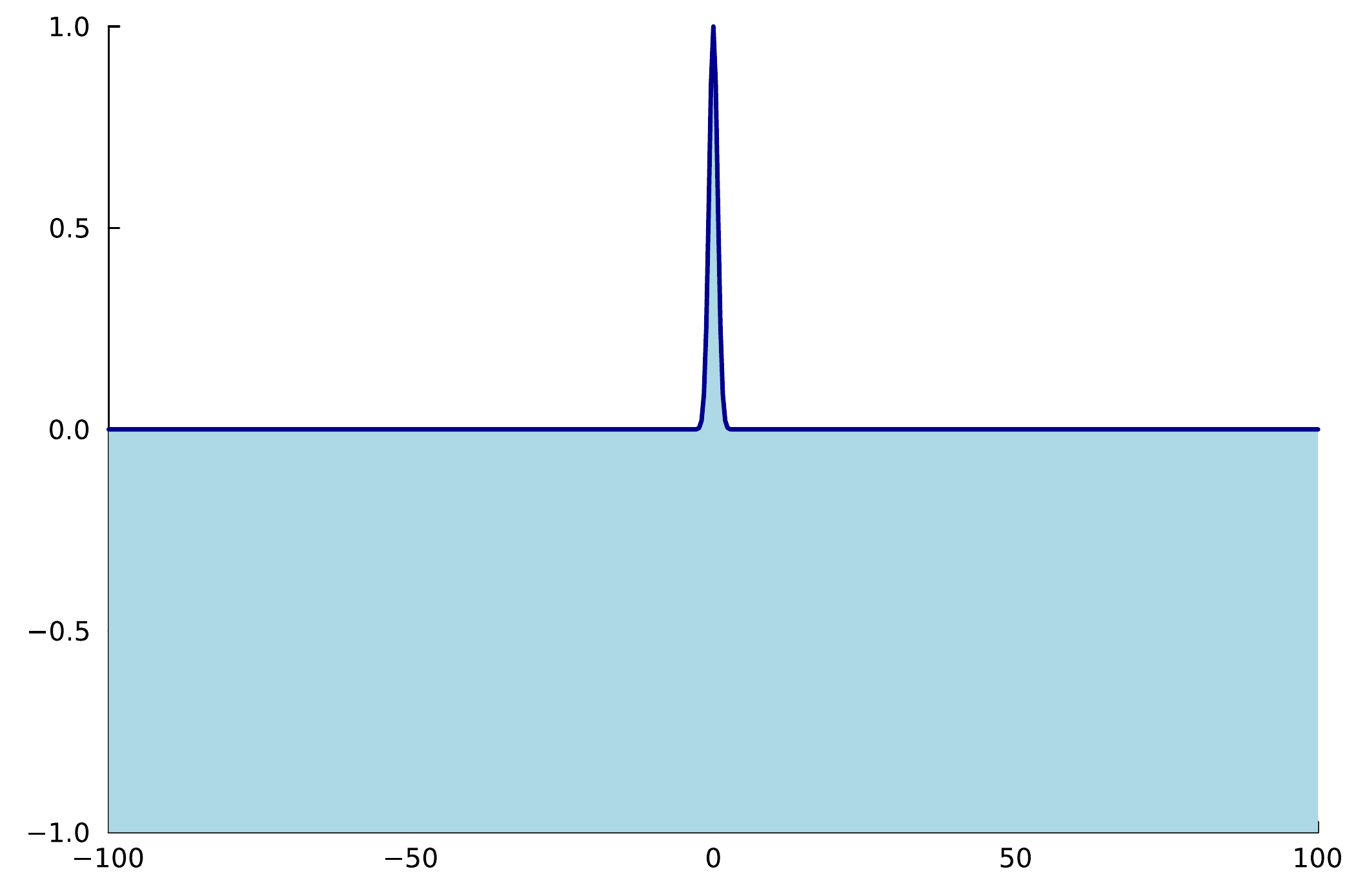}
  \caption{Initial condition, \(q(x,0) = \ee^{-x^2}\).}
  \label{fig:sfig185}
\end{subfigure}%
\begin{subfigure}{.5\textwidth}
  \centering
  \includegraphics[width=.7\linewidth]{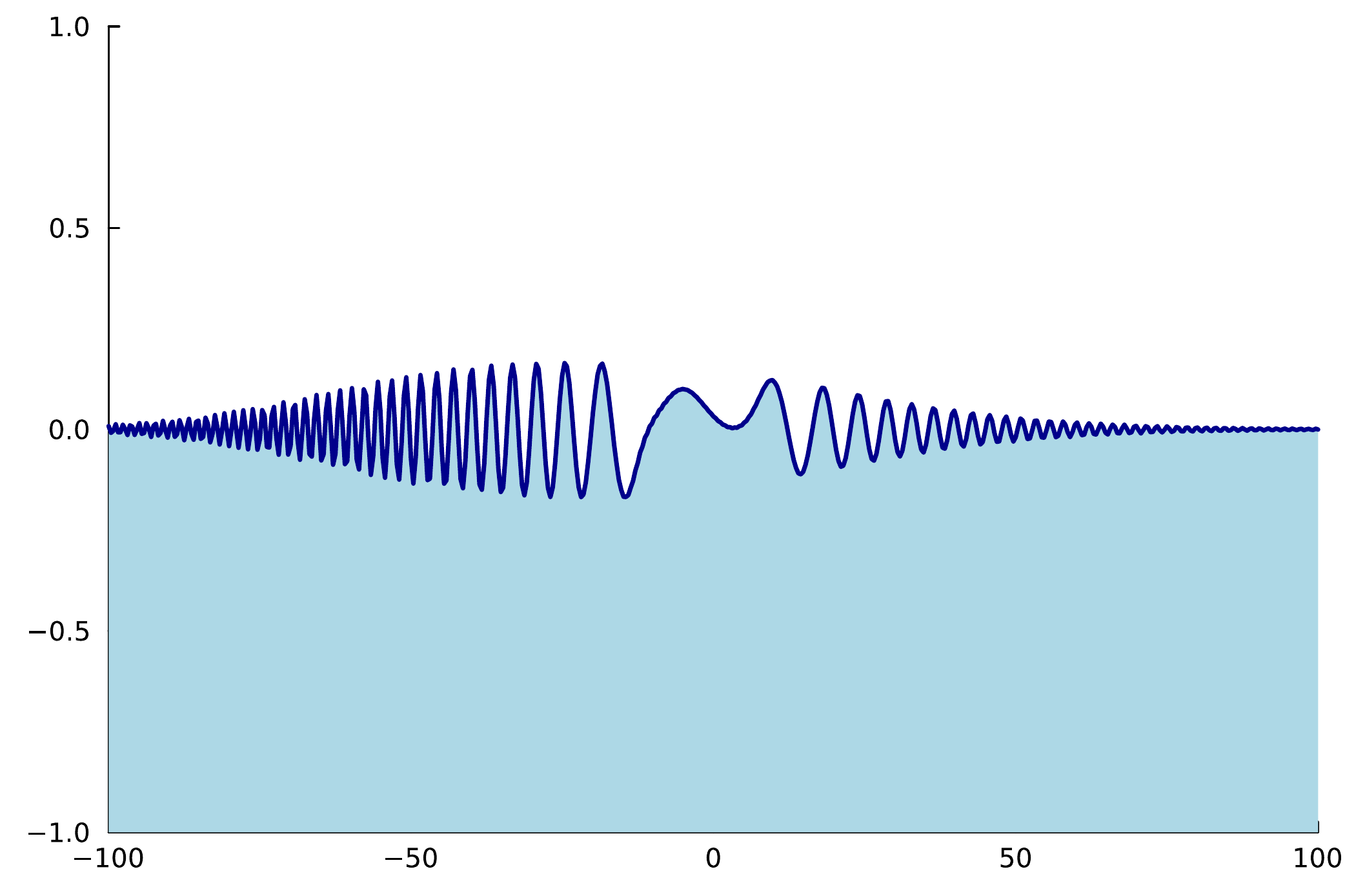}
  \caption{Solution at \(t = 10\).}
  \label{fig:sfig270}
\end{subfigure}
\caption{Solution to the Eckhaus equation \cite{eckhaus}: \(\ii q_t+q_{xx}+2|q|^2_xq+|q|^4q=0\). \textbf{Parameters:} \(L = 200\), \(m = 2^{10}\), \(\triangle t = 0.01\), \(\sigma(x) = \) \eqref{eq:sigma}, \(\gamma(x) = \) \eqref{eq:newgamma}, \(k_1 = 0\), \(f_2 = 1000\). \textbf{Runtime:} 5 seconds.}
\label{fig:eckhaus}
\end{figure}

In this section, we include plots of solutions to different problems computed using our damping method. We begin by solving the KdV equation \eqref{eq:kdv} with an initial condition that produces a two-soliton solution, shown in \eqref{fig:twosoliton}. In \eqref{fig:riemannappend}, we solve the Riemann problem for the KdV equation \eqref{eq:cauchy} with an initial condition motivated by an example in \cite{cauchyproblem}. Qualitative comparison of the solution to that of Grava and Klein confirms its accuracy. We then solve the Riemann problem for the KdV equation with a soliton on top of a shelf for the initial condition in \eqref{fig:solitonshelf}. In \eqref{fig:eckhaus}, we plot the real part of the solution to the Eckhaus equation \cite{eckhaus}. The Eckhaus equation models wave propagation in dispersive media and falls within the nonlinear Schrödinger class of equations.\\

\clearpage 

\bibliographystyle{amsplain}
\bibliography{references}

\end{document}